\documentclass{conm-p-l}

\newtheorem{theorem}{Theorem}[section]
\newtheorem{lemma}[theorem]{Lemma}
\newtheorem{proposition}[theorem]{Proposition}
\newtheorem{corollary}[theorem]{Corollary}

\theoremstyle{definition}
\newtheorem{definition}[theorem]{Definition}
\newtheorem{example}[theorem]{Example}

\theoremstyle{remark}
\newtheorem{remark}[theorem]{Remark}

\numberwithin{equation}{section}

\def \End{{\rm End}}

\def \im {{\rm Im}}
\def \deg {{\rm deg}}
\def \g {{\mathfrak g}}

\def \im {{\rm im}\;}

\def \Res{{\rm Res}}
\def \wt {{\rm wt}}
\def \Z{{\mathbb Z}}
\def \N{{\mathbb N}}
\def \C{{\mathbb C}}
\def \l{{\lambda}}

\def \bex{\begin{example}\label}
\def \eex{\end{example}}
\def \bl{\begin{lemma}\label}
\def \el{\end{lemma}}
\def \bc{\begin{corollary}\label}
\def \ec{\end{corollary}}
\def \bd{\begin{definition}\label}
\def \ed{\end{definition}}
\def \bt{\begin{theorem}\label}
\def \et{\end{theorem}}
\def \bp{\begin{proposition}\label}
\def \ep{\end{proposition}}
\def \br{\begin{remark}\label}
\def \er{\end{remark}}

\def \<{\langle}
\def \>{\rangle}

\begin{document}

\title{Vertex Lie algebras, vertex Poisson algebras and
vertex algebras}

\author{Chongying Dong}
\address{Department of Mathematics, University of California,
Santa Cruz, CA 95064}
\email{dong@math.ucsc.edu}
\thanks{C. Dong was supported by NSF grant DMS-9303374, DMS-9987656 and a
research grant from the Committee on Research, UC Santa Cruz.
H. Li was supported by NSF grant DMS-9616630.
G. Mason was supported by NSF grant 
DMS-9401272, DMS-9700909 and a research grant from the Committee on Research, UC
Santa Cruz.}

\author{Haisheng Li}
\address{Department of Mathematical Sciences, 
Rutgers University, Camden, NJ 08102}
\email{hli@crab.rutgers.edu}

\author{Geoffrey Mason}
\address{Department of Mathematics, University of California,
Santa Cruz, CA 95064}
\email{gem@math.ucsc.edu}

\subjclass{Primary 17B69; Secondary 16A68, 81T40}


\begin{abstract}
The notions of vertex Lie algebra and vertex Poisson algebra
are presented  and connections among 
vertex Lie algebras, vertex Poisson algebras and vertex algebras are
discussed.
\end{abstract}

\maketitle

\section{Introduction}

Vertex (operator) algebras (see [B], [FLM]) which 
are a family of new ``algebras'' are known essentially to be chiral algebras 
in two-dimensional quantum field theory. The new algebras are
closely related to classical Lie algebras of certain types.  
On one hand, from [B] one can obtain a Lie algebra $V/DV$ from
any vertex algebra $V$ where $D$ is a canonical endomorphism of $V$. 
(The famous Monster Lie algebra
is a subalgebra of $V/DV$ for a suitable $V.$) On the other hand, one
can construct vertex (operator) algebras from certain highest weight
representations of familiar infinite dimensional Lie algebras such as
the Virasoro algebra, affine Kac-Moody algebras, Heisenberg algebras
(cf. [DL], [FF1], [FZ], [Li2], [MP]).

In this paper we define a notion of what we call vertex Lie algebra
to unify the familiar infinite dimensional Lie algebras.
The definition of a vertex Lie algebra is motivated by 
the Virasoro algebra, affine algebras and the notion of vertex algebra. 
Roughly speaking, the notion of vertex Lie algebra is a ``stringy'' 
analogue of the notion of Lie algebra, which generalizes
the Virasoro algebra and affine Lie algebras. 
A vertex Lie algebra defined  in this paper
is a Lie algebra whose underlying vector space is essentially 
the loop space of a certain vector space $U$
and it is a kind of ``affine algebra''
based on vector space $U$ instead of a finite dimensional Lie algebra.
So for every vector $u$ in the base vector space, one can form
the field or vertex operator $u(z)=\sum_{n\in\Z}u(n)z^{-n-1}.$ 
In the case of the Virasoro algebra, $U$ is either one dimensional
or two dimensional depending on the center being zero or not.
(The Virasoro algebra is not a classical affine algebra in any sense.)
We also  define a notion of 
vertex Poisson algebra. We study the connection among
vertex Lie algebras, vertex Poisson algebras and vertex algebras.

The notion of vertex Lie algebras not only unifies the Virasoro
algebra, affine Lie algebras, loop algebras and other important
infinite dimensional Lie algebras, but also provides new examples
of vertex algebras via their ``highest weight representations.''
We prove that for each complex number there is a 
``highest weight module'' (based a polar decomposition) 
for the vertex Lie algebra which has the structure of a vertex algebra.
We also show how a restricted module for the vertex Lie algebra
has a natural module structure for certain vertex algebras
constructed from the vertex Lie algebra. So it is very important 
to construct new vertex Lie algebras.

Vertex Poisson Lie algebras are a special class of vertex Lie algebras
whose base vector spaces are  Poisson algebras satisfying additional
axioms. Closely related to vertex Poisson algebras are
Possion brackets studied in [DFN], [DN] and [GD]. 
We were amazed to find that
local Poisson brackets introduced in [DFN] was so close to Borcherds'
commutator formula in the theory of vertex algebras.
Our first exercise is to make local Poisson brackets
precise in terms of formal variables where delta functions are
formal series instead of distributions.
Differential geometric Poisson brackets introduced in [DN], where 
interesting connections between differential 
geometric objects (connection, curvature and symplectic structure)
and algebraic objects (Lie algebra, commutative associative algebra)
have been established, provide a lot of examples of vertex Poisson
Lie algebras. At the end of Section 3, we quote several interesting 
results from [DFN], [DN] and [Po] regarding to differential 
geometric Poisson brackets.

In [BD], a notion of cossion algebra was defined in terms of 
algebraic geometry where as designated by the authors, 
cossion  is the combination of the two words chiral and poission. 
Later, a notion of vertex Poisson algebra was defined in [EF]
in terms of associative rings. (The relation between coission algebras and 
vertex Poisson algebras was discussed in [EF].) The notion of 
vertex Poisson algebra presented  here is defined in terms of formal calculus
and classical (Lie and associative) algebras. Presumably, our notion is 
essentially the same as that of [EF].

Most of this work was carried out in the late 96 and early 97. 
During this time, Kac's book [K2] appeared where
similar results to ours had been obtained. In particular, 
a notion of conformal algebra was introduced. 
A conformal algebra satisfies a set of axioms which are certain 
modifications of those for a vertex algebra and a conformal algebra
naturally gives to a vertex Lie algebra. A notion of vertex Lie algebra 
was also introduced in [Pr]. It seems that
the notion of vertex Lie algebra in [Pr] the notion 
of conformal algebra [K2] are the same. There are certain overlaps between 
this work and [Pr]. In particular, Theorem 4.8 and Lemma 5.3 were also 
obtained in [Pr].

The paper is organized as follows. Section 2 is about formal calculus.
In Section 3 we formulate the notions of local vertex Lie algebra
and local vertex Poisson algebra and present some results. 
We also give several examples including
the Witt algebra, Virasoro algebra, loop algebras, affine algebras
to illustrate the concepts. In particular, we discuss vertex Poisson algebras
for the base vector spaces being the Poisson algebras
$S(\g)$ which is the symmetric algebra of Lie algebra $\g.$ 
In Section 4 we construct vertex algebras and their modules
from vertex Lie algebras by using the frame work of local system introduced 
in [Li2]. 
Since it has been an important problem in the theory
of vertex algebras to construct new examples of vertex algebras, 
the vertex algebra construction based on vertex Lie algebras
is of certain importance. 
Section 5 is devoted to the study of Lie algebras
and Poisson algebras associated to vertex algebras. Using the results
obtained in Section 4 we construct a vertex algebra $V^{[l]}$ 
for any vertex algebra $V$ and any complex number $l.$ In particular,
if $V$ is a vertex operator algebra of CFT type we construct 
a new vertex operator algebra $V^{[l]}$ whose central charge is
$lc$ where $c$ is the central charge of $V.$ We also investigate
the Poisson algebra $P_2(V)$ (see [Zh]) for vertex algebra constructed
from a vertex Lie algebra. We show that for any Lie algebra
$\g$ there exists a vertex algebra $V$ such that $P_2(V)$ is
isomorphic to $S(\g)$ as Poisson algebras.

\section{Calculus of formal variables}

In this section we review some formal variable notations and
the fundamental properties of delta function from [FLM] and
we also formulate some simple results which will be used later.  

Throughout the whole paper, $x, x_{0},x_{1},\ldots, y, y_{0},y_{1},
\ldots, z,z_{0},z_{1},\ldots$ are independent commuting formal variables.
The symbols $\Z, \N, \C$ stand for the integers, the nonnegative 
integers and the complex numbers. Vector spaces are over the ground 
field $\C$.

For any vector space $U$, following [FLM] we set
\begin{eqnarray}
U[[x,x^{-1}]]&=&\left\{ \sum_{n\in {\Z}}a_{n}x^{n}\;\;|\;\;a_{n}\in U\right\},
\\
U((x))&=&\left\{\sum_{n\ge m}a_{n}x^{n}\;\;|\;\;m\in {\Z},\; a_{n}
\in U\right\},\\
U[[x]]&=&\left\{\sum_{n=0}^{\infty}a_{n}x^{n}\;\;|\;\;a_{n}\in U\right\},\\
U[x,x^{-1}]&=&\left\{\sum_{n=m}^{k}a_{n}x^{n}\;\;|\;m,k\in\Z,\;a_{n}
\in U\right\}.
\end{eqnarray}

The formal delta-function is defined to be the formal series 
$\delta(x)=\sum_{n\in {\Z}}x^{n}$ (an element of ${\C}[[x,x^{-1}]]$)
 and its fundamental property is
\begin{eqnarray}
f(x)\delta(x)=f(1)\delta(x)\;\;\mbox{ for }f(x)\in U[x,x^{-1}].
\end{eqnarray}
Furthermore
\begin{eqnarray}
f(x)\delta\left(\frac{x}{y}\right)=f(y)\delta\left(\frac{x}{y}\right)\;\;\;
\mbox{ for }f(x)\in U[[x,x^{-1}]].
\end{eqnarray}
For any $n\in {\Z}$, $(x+y)^{n}$ is defined to be the formal series 
$\sum_{i=0}^{\infty}\binom{n}{i}x^{n-i}y^{i}$, 
where $\binom{n}{i}=\frac{1}{i!}n(n-1)\cdots (n+1-i)$.
We also have the following identities in three formal variables (see [FLM]):
\begin{eqnarray}
& &x^{-1}\delta\left(\frac{y-z}{x}\right)
=y^{-1}\delta\left(\frac{x+z}{y}\right),\label{edelta1}\\
& &x^{-1}\delta\left(\frac{y-z}{x}\right)
- -x^{-1}\delta\left(\frac{z-y}{-x}\right)
=z^{-1}\delta\left(\frac{y-x}{z}\right).
\end{eqnarray}
By using the well-known Taylor formula and (\ref{edelta1}) we get
$$e^{-z\partial/\partial y}x^{-1}\delta\left(\frac{y}{x}\right)=
x^{-1}\delta\left(\frac{y-z}{x}\right)=y^{-1}\delta\left(\frac{x+z}{y}\right)
=e^{z\partial/\partial x}y^{-1}\delta\left(\frac{x}{y}\right).$$
This amounts to that
\begin{eqnarray}\label{epart}
\left(-\dfrac{\partial}{\partial x}\right)^{k}
y^{-1}\delta\left(\frac{x}{y}\right)
=\left(\dfrac{\partial}{\partial y}\right)^{k}
y^{-1}\delta\left(\frac{x}{y}\right)\;\;\;\mbox{ for }k\in \N.
\end{eqnarray}
Later, we shall frequently use the derivatives of 
$y^{-1}\delta(\frac{x}{y})$. 
For convenience, we set
\begin{eqnarray}
\Delta^{(k)}(x,y)=\left(\dfrac{\partial}{\partial x}\right)^{k}
y^{-1}\delta\left(\frac{x}{y}\right)
=\left(-\dfrac{\partial}{\partial y}\right)^{k}y^{-1}\delta\left(\frac{x}{y}
\right)
\end{eqnarray}
for $k\in \N.$ If $k=0$ we simply write $\Delta(x,y)$ for $\Delta^{(0)}(x,y)
=y^{-1}\delta(\frac{x}{y}).$

For $f(x) \in \C[[x,x^{-1}]]$, since $f(x)\Delta(x,y)=f(y)\Delta(x,y)$, 
we have
\begin{eqnarray}\label{echange}
f(x)\Delta^{(k)}(x,y)&=&\left(-\dfrac{\partial}{\partial y}\right)^{k}
f(x)\Delta(x,y)\nonumber\\
&=&(-1)^{k}\left(\dfrac{\partial}{\partial y}\right)^{k}
\left(f(y)\Delta(x,y)\right)
\nonumber\\
&=&(-1)^{k}\sum_{j=0}^{k}\binom{k}{j}f^{(k-j)}(y)\Delta^{(j)}(x,y).
\end{eqnarray}

\bl{lformula0}
 Let $m,n\in \N$. If $m\le n$ we have
\begin{eqnarray}\label{eformula1}
(x-y)^{m}\Delta^{(n)}(x,y)=\binom{-n}{m}m!\Delta^{(n-m)}(x,y).
\end{eqnarray}
If $m>n$ we have
\begin{eqnarray}\label{eformula2}
(x-y)^{m}\Delta^{(n)}(x,y)=0.
\end{eqnarray}
\el

\begin{proof}
Since $(x-y)\Delta(x,y)=0$, we have
$$0=\left(\dfrac{\partial}{\partial x}\right)^{n}\left((x-y)\Delta(x,y)\right)
=(x-y)\Delta^{(n)}(x,y)+n\Delta^{(n-1)}(x,y).$$
Thus $(x-y)\Delta^{(n)}(x,y)=-n\Delta^{(n-1)}(x,y)$ for $n\ge 1$.
Then (\ref{eformula1}) follows from induction on $m$ immediately.
Furthermore, (\ref{eformula2}) follows from (\ref{eformula1}) and the fact
that $(x-y)\Delta^{(0)}(x,y)=0.$
\end{proof}

Recall the following result from [Li2] (cf. [FLM], Proposition 8.1.3):

\bl{lbasis}
Let $U$ be any vector space and let $f_{i}(y)\in U[[y,y^{-1}]]$ for 
$i=0,1,\ldots, n$. Then 
\begin{eqnarray}
f_{0}(y)\Delta(x,y)+f_{1}(y)\Delta^{(1)}(x,y)+\cdots +f_{n}(y)\Delta^{(n)}(x,y)
=0
\end{eqnarray}
if and only if $f_{i}(y)=0$ for all $i$. Consequently, the expression of an 
element $h(x,y)$ of $U[[x,y,x^{-1},y^{-1}]]$ as a finite sum 
$\sum_{i=0}^{n}g_{i}(y)\Delta^{(i)}(x,y)$ is unique if it exists.
\el

\br{rbasis}
Similarly, if $U$ is a vector space and $f_{i}(x)\in U[[x,x^{-1}]], 
i=0,1,\ldots n$, then
\begin{eqnarray}
f_{0}(x)\Delta(x,y)+f_{1}(x)\Delta^{(1)}(x,y)+\cdots +f_{n}(x)\Delta^{(n)}(x,y)
=0
\end{eqnarray}
if and only if $f_{i}(x)=0$ for all $i$. Consequently, the expression of an 
element $h(x,y)$ of $U[[x,y,x^{-1},y^{-1}]]$ as a finite sum 
$\sum_{i=0}^{n}g_{i}(x)\Delta^{(i)}(x,y)$ is unique if it exists.
\er

The following result can also be found in [K2] (cf. [FLM], Proposition 8.1.3): 
(Our proof is slightly different from that of [K2].)

\bp{plocal}
Let $U$ be any vector space and let $f(x,y)\!\in\! U[[x,y,x^{-1},y^{-1}]]$, 
a formal series in $x,y$ with coefficients in $U$. Then $(x-y)^{k+1}f(x,y)=0$ 
for some nonnegative integer $k$ if and only if 
$\displaystyle{f(x,y)=\sum_{i=0}^{k}f_{i}(y)\Delta^{(i)}(x,y)}$
for some $f_{0}(y),\ldots, f_{k}(y)\in U[[y,y^{-1}]]$.
\ep

\begin{proof}
Since $(x-y)^{m}\Delta^{(n)}(x,y)=0$ for any nonnegative 
integers $m>n$, one direction is clear.
We shall prove the other direction by using induction on $k$.
Suppose that $(x-y)f(x,y)=0$. Let $f(x,y)=\sum_{m,n\in {\Z}}a(m,n)x^{m}y^{n}$.
Then we have $a(m+1,n)=a(m,n+1)$ for any $m,n\in {\Z}$, so that
$a(m,n)=a(m+n,0)$ for $m,n\in {\Z}$. Then
$f(x,y)=f_{0}(y)\delta\left(\frac{x}{y}\right)$, where 
$f_{0}(y)=\sum_{n\in {\Z}}a(n,0)y^{n}$. Thus it is true for $k=0$.

Now, suppose  that it is true for $k$ and that a formal series $f(x,y)$ 
satisfies the relation $(x-y)^{k+2}f(x,y)=0$.
Set $g(x,y)=(x-y)f(x,y)$. Then $(x-y)^{k+1}g(x,y)=0$. By the inductive 
hypothesis, there are $g_{0}(y),\ldots, g_{k}(y)\in U[[y,y^{-1}]]$ such that
$$g(x,y)=\sum_{j=0}^{k}g_{i}(y)\Delta^{(j)}(x,y).$$
Set 
$$F(x,y)=\sum_{j=0}^{k}\frac{1}{j+1}g_{j}(y)
\Delta^{(j+1)}(x,y).$$ 
Since $(x-y)\Delta^{(n)}(x,y)=-n\Delta^{(n-1)}(x,y)$
for any positive integer $n$, we have
\begin{eqnarray*}
& &(x-y)(f(x,y)+F(x,y))\\
&=&\sum_{j=0}^{k}\left( g_{j}(y)\Delta^{(j)}(x,y)+
\frac{1}{j+1}(x-y)g_{j}(y)\Delta^{(j+1)}(x,y)\right)=0.
\end{eqnarray*}
{}From the base step, we have $f(x,y)+F(x,y)=g(y)\Delta(x,y)$ 
for some $g(y)\in U[[y,y^{-1}]]$. Then the inductive step follows 
immediately. This concludes the proof.
\end{proof}

\br{rcoprodbasis}
Set
\begin{eqnarray}
A(x,y)=\sum_{n\in \N}\C[[x,x^{-1}]]\Delta^{(n)}(x,y)
\subset 
\C[[x,y,x^{-1},y^{-1}]].
\end{eqnarray}
{}From (\ref{echange}) we also have
\begin{eqnarray}
A(x,y)=\sum_{n\in \N}\C[[y,y^{-1}]]\Delta^{(n)}(x,y).
\end{eqnarray}
In view of Proposition \ref{plocal}, $A(x,y)$ can be canonically 
defined by
$$A(x,y)=\{ f(x,y)\in \C[[x,y,x^{-1},y^{-1}]]\;|\; (x-y)^{k}f(x,y)=0
\;\;\;\mbox{ for some }k\in \N\}.$$
It follows from Lemma \ref{lbasis} (with $U=\C$) that
\begin{eqnarray}
A(x,y)=\bigoplus_{n\in \N}\C[[x,x^{-1}]]\Delta^{(n)}(x,y)
=\bigoplus_{n\in \N}\C[[y,y^{-1}]]\Delta^{(n)}(x,y).
\end{eqnarray}
\er

We shall need the following result (cf. [MN], Lemma 1.1.1):

\bl{lele}
Let $U$ be a vector space and let 
$$f(x,y)\in \C((x))\otimes U[[y,y^{-1}]],\;\mbox{ or }
\C((y))\otimes U[[x,x^{-1}]].$$
Then $(x-y)^{k}f(x,y)=0$ for some nonnegative integer $k$
if and only if $f(x,y)=0$.
\el

\begin{proof}
Suppose that $f(x,y)\in \C((x))\otimes U[[y,y^{-1}]]$.
Then, for {\em any} $n\in \Z$,
$$(-y+x)^{n}f(x,y)\;\;\;\mbox{ exists }$$
in $U[[x,x^{-1},y,y^{-1}]]$. Then from [FLM] (Chapter 2) we have
\begin{eqnarray}
(-y+x)^{m}\left((-y+x)^{n}f(x,y)\right)=(-y+x)^{m+n}f(x,y)
\end{eqnarray}
for any $m,n\in \Z$. It follows immediately that 
$(x-y)^{k}f(x,y)=0$ for some nonnegative integer $k$
if and only if $f(x,y)=0$. Similarly, this is true for
$f(x,y)\in \C((y))\otimes U[[x,x^{-1}]]$.
\end{proof}

\br{rlocal}
Let $L$ be any Lie algebra and let $\psi(x), \phi(x)\in L[[x,x^{-1}]]$.
Then by Proposition \ref{plocal}, there exists a nonnegative integer 
$k$ such that
\begin{eqnarray}\label{elocal1}
(x-y)^{k}\psi(x)\phi(y)=(x-y)^{k}\phi(y)\psi(x)
\end{eqnarray}
 if and only if
there exist $f_{0}(y),\ldots, f_{k}(y)\in L[[y,y^{-1}]]$ such that
\begin{eqnarray}\label{egl}
[\psi(x),\phi(y)]=\sum_{i=0}^{k}f_{i}(y)\Delta^{(i)}(x,y).
\end{eqnarray}
Furthermore, it follows from Lemma \ref{lbasis} that such 
$f_{0}(y),\ldots, f_{k}(y)$ are uniquely determined by $\psi(x)$ and $\phi(x)$.
Suppose $\psi(x)\in L[x,x^{-1}]$. Then it follows from Lemma \ref{lele}
that (\ref{elocal1}) holds if and only if $\psi(x)\phi(y)=\phi(y)\psi(x)$.
\er

\section{Local vertex Lie algebra and vertex Poisson algebras}

In this section we formulate certain notions of local vertex 
Lie algebra and  local vertex Poisson algebra, and give
several consequences of the definition. 
A notion of a vertex Lie algebra was also independently defined in [Pr]
and [EF], and a notion of a vertex Poisson algebra has been defined 
in [EF] (and in [BD]). The notion of a conformal algebra defined in [K2]
is also closely related to the notion of a vertex Lie algebra.

\bd{dvla} (1)
A {\em local vertex Lie algebra } is a quadruple $(L,U,d,\rho)$ 
consisting of a
Lie algebra $L$, a vector space $U$ and a linear map $\rho$ from
$L(U)=U\otimes {\C}[t,t^{-1}]$ by definition onto $L$ such that
$\ker \rho= \im\; \hat{d}$, where $\hat{d}=d\otimes 1+1\otimes \dfrac{d}{dt}$
and $d$ is a partially defined linear map from $U$ to $U$, and that
for any $u,v\in U$, there exist finitely many 
$f_{0}(u,v),\ldots, f_{r}(u,v)\in U$ and nonnegative integers 
$k_{0},\ldots, k_{r},$ $l_{0},\ldots, l_{r}$
(where $r$ depends on both $u$ and $v$) such that
\begin{eqnarray}\label{edlva}
\ \ \ \ \ [u(x), v(y)]=f_{0}(u,v)^{(k_{0})}(y)\Delta^{(l_0)}(x,y)+
\cdots +
f_{r}(u,v)^{(k_{r})}(y)\Delta^{(l_r)}(x,y),
\end{eqnarray}
where for $i\in \N,$ $u\in U,$ 
$$u^{(i)}(x)=\left(\dfrac{d}{dx}\right)^{i}\sum_{n\in {\Z}}u(n)x^{-n-1}=
\sum_{n\in {\Z}}\binom{-n-1}{i}i!u(n)x^{-n-i-1}$$ and 
$u(n)=\rho (u\otimes t^{n}).$ We simply write
$L=(L,U,d,\rho).$

(2) A ${\Z}$-{\em graded} vertex Lie algebra is a vertex Lie algebra 
$L$ over $U$ such that $L$ is a ${\Z}$-graded Lie algebra and
$U=\oplus_{n\in {\Z}}U_{(n)}$ is a ${\Z}$-graded space such that
$\deg\; u(n)=m-n-1$ for $u\in U_{(m)},\; n\in {\Z}$.
\ed

We shall conventionally call $L$ {\em a local vertex Lie algebra 
over the base space $U$}. 

In the following example, we shall see that affine Lie algebras
and the Virasoro algebra are local vertex Lie algebras. 

\bex{exa1} (1) Let $Witt=\oplus_{n\in {\Z}}{\C}L(n)$ be the Witt Lie algebra
or the centerless Virasoro algebra, where
\begin{eqnarray}
[L(m),L(n)]=(m-n)L(m+n)\;\;\;\mbox{ for }m,n\in \Z.
\end{eqnarray}
Then $Witt$ is a local vertex Lie
algebra over the base space $U={\C}\omega$ where $\rho(\omega\otimes t^{m})
=\omega(m+1)=L(m)$ 
for $m\in {\Z}$, the domain of $d$ is $0$ and
\begin{eqnarray}
[\omega(x),\omega(y)]=\omega^{(1)}(y)\Delta(x,y)-2\omega(y)\Delta^{(1)}(x,y).
\end{eqnarray}

(2)  The Virasoro algebra $Vir=\oplus_{n\in {\Z}}{\C}L(n)\oplus \C c$ 
is a local vertex Lie algebra over the base space 
$U={\C}\omega\oplus {\C}c$, where
$$\rho(\omega\otimes t^{m})=\omega(m+1)
=L(m),\;\; \rho(c\otimes t^{m})=c(m)=\delta_{m,-1}c$$
for $m\in {\Z}$, the domain of $d$ is ${\C}c$ with $d=0$
and
\begin{eqnarray}
& &[\omega(x),\omega(y)]
=\omega^{(1)}(y)\Delta(x,y)-2\omega(y)\Delta^{(1)}(x,y)
- -\frac{1}{12}c(y)\Delta^{(3)}(x,y),\\
& &[\omega (x), c(y)]=0.
\end{eqnarray}

(3) Let $\g$ be a Lie algebra and let $L(\g)$ be the 
loop (Lie) algebra
\begin{eqnarray}
L(\g)=\g\otimes {\C}[t,t^{-1}],
\end{eqnarray}
where
\begin{eqnarray}
[a\otimes t^{m},b\otimes t^{n}]=[a,b]\otimes t^{m+n}
\;\;\;\mbox{ for }a,b\in \g,\; m,n\in \Z.
\end{eqnarray} 
The loop algebra $L(\g)$ is a local vertex Lie algebra over the 
base space $\g$ with $\rho$ being the identity map and the domain of $d$ 
being zero, where
\begin{eqnarray}
[a(x),b(y)]=[a,b](y)\Delta(x,y)\;\;\;\mbox{ for }a,b\in \g.
\end{eqnarray}

(4) If $\g$ is a Lie algebra  with a symmetric nondegenerate invariant 
bilinear form $(\cdot|\cdot)$, then the corresponding
affine Lie algebra $\hat{\g}=L(\g)\oplus {\C}c$ is a local vertex Lie algebra
over the base space $U=\g\oplus {\C}c$ where
$$\rho(a\otimes t^{m})
=a\otimes t^{m},\;\; \rho(c\otimes t^{m})=\delta_{m,-1}c\;\;\;\mbox{ for }
a\in \g,\; m\in {\Z},$$
the domain of $d$ is ${\C}c$ with $d=0,$  and
\begin{eqnarray}
& &[a(x),b(y)]=[a,b](y)\Delta(x,y)-(a|b)c(y)\Delta^{(1)}(x,y),\\
& &[a(x), c(y)]=0\;\;\;\mbox{ for }a,b\in \g.
\end{eqnarray}
It is not hard to see that all the local vertex Lie algebras discussed here
are graded. 
\eex

\br{rdvla}
(1) Let $(L,U,d,\rho)$ be a local vertex Lie algebra.
Then $L\simeq L(U)/\ker \rho$. Clearly,  
$u\otimes t^{-1}\notin \im d=\ker\rho$ for any $0\ne u\in U.$ Thus $u(-1)\ne 0$
and $u(x)\ne 0$. 
Then the linear map $d$ is uniquely determined by the condition 
$\dfrac{d}{dx}a(x)=(da)(x)$.
 Set $U^{0}=\ker d$. 
Then $U^{0}\otimes \frac{d}{dt}{\C}[t,t^{-1}]\subset \im d=\ker \rho$.

(2) The operator $d$ implies the following:  
In the definition of local vertex Lie algebra, if $f_i(u,v), df_i(u,v),
\ldots, d^{k_i}f_i(u,v)$ are defined then $f_i(u,v)^{(k_i)}(y)$ 
can be replaced by $(d^{k_i}f_i(u,v))(y).$ In particular,
if $d$ is defined on whole $U,$ we can choose $f_i(u,v)$ so that
$k_i=0$ for all $i.$
\er

The following are consequences of the definitions.

\bp{central} (1) Let $(L,U,d,\rho)$ be a local vertex Lie algebra with 
$\ker d=U^{0}$. Then
\begin{eqnarray}
u(n)=0\;\;\;\mbox{ for }u\in U^{0},\; n\ne -1,
\end{eqnarray}
and $u(-1)$ for $u\in U^{0}$ form a central subalgebra of $L$.

(2) The components of $u(m)$ and $v(n)$ of $u(x_1)$ and $v(x_2)$
have bracket
\begin{equation}\label{ecomcom}
[u(m),v(n)]\!=\!\sum_{i=0}^{r}\binom{m}{l_i}
\binom{m\!+\!n\!-\!l_i}{k_{i}}(-1)^{l_i+k_i}l_i!k_{i}!
f_{i}(u,v)(\!m\!+\!n-\!l_i\!-\!k_{i}\!).
\end{equation}

(3) Set $L^{-}=\rho(U\otimes t^{-1}{\C}[t^{-1}])$,
$L^{+}=\rho(U\otimes {\C}[t])$ and $L^{0}=\rho(U)$ (a subspace 
of $L^{+}$). Then $L^{\pm}$ and $L^{0}$ are 
Lie subalgebras of $L$ and $L=L^{+}\oplus L^{-}$ is a polar decomposition.
Furthermore, $L^{0}=U/(U^{0}+\im d)$ and 
$\rho$ is a linear isomorphism from 
$(U^{0})'\otimes t^{-1}\oplus U'\otimes {\C}[t,t^{-1}]$ onto
$L$, where $U'$ is a subspace of $U$ such that $U=(U^{0}+\im d)\oplus U'$
and $(U^{0})'$ is a subspace of $U^{0}$ such that $(U^{0})'\cap \im d=0$.
If $\ker d \cap \im d=0$, then we may take $(U^{0})'=U^{0}$.

(4) If $L=\oplus_{n\in {\Z}}L_{(n)}$ is a graded local vertex Lie algebra, 
we have a triangular 
decomposition: $L=L_{+}\oplus L_{0}\oplus L_{-}$, where
$L_{\pm}=\oplus_{n=1}^{\infty}L_{(\pm n)}$ and $L_{0}=L_{(0)}$.
\ep

\begin{proof} Since $\ker \rho =\im \hat{d}$, for $u\in U$,
$u=u\otimes t^{0}\in \ker \rho$ if and only if
\begin{eqnarray}\label{eexplicit1}
u\otimes t^{0}=dv\otimes t^{n}+nv\otimes t^{n-1}
\end{eqnarray}
for some $v\in U,\; n\in \Z$. Clearly,  (\ref{eexplicit1})
is equivalent to that either $u\otimes t^{0}=dv\otimes t^{n},\;\; n=0$, or
$dv=0$ and $u\otimes t^{0}=nv\otimes t^{n-1}$.
Then $u=u\otimes t^{0}\in \ker \rho$ if and only if $u\in U^{0}+\im d$.
Thus  $L^{0}=U/(U^{0}+\im d)$.
>From Remark \ref{rdvla} (1), we have
\begin{eqnarray}\label{edum=-mum-1}
(du)(m)=-mu(m-1)\;\;\;\mbox{ for }u\in U,\; m\in \Z,
\end{eqnarray}
assuming that $u$ is in the domain of $d$.
Then, for $u\in U^0$, $u(m)=0$ for $m\ne -1$, i.e., $u(x)=u(-1)\in L[x,x^{-1}]$.
It follows from Remark \ref{rlocal} that $u(-1)$ commutes with $v(x)$ for
any $v\in U.$ This proves (1). (2) is immediate from 
(\ref{edlva}) by considering the coefficients of $x^{-m-1}y^{-n-1}.$
The first part of (3) follows from (2) by noticing that 
$\binom{m}{l_i}\binom{m+n-l_i}{k_i}=0$
if $m+n-l_i-k_i<0$ and $m,n\ge 0$. 
{}From (\ref{edum=-mum-1}), we have
$$\rho(U^{0}\otimes t^{-1}+U'\otimes {\C}[t,t^{-1}])=L.$$
{}From a simple fact in linear algebra, it suffices to  prove 
$$\ker \rho \cap (U'\otimes {\C}[t,t^{-1}])=0,$$
that is,
$$\im \hat{d} \cap (U'\otimes {\C}[t,t^{-1}])=0.$$
Otherwise, there exist $0\ne u^{i}\in U',\; m_{1}>\cdots >m_{r}$ such that
$$ u^{1}\otimes t^{m_{1}}+\cdots +u^{r}\otimes t^{m_{r}}\in \im \hat{d}.$$
Then
\begin{eqnarray}\label{ebasiscalculation}
u^{1}\otimes t^{m_{1}}+\cdots +u^{r}\otimes t^{m_{r}}
=\sum_{j=1}^{s}( (dv^{j}\otimes t^{n_{j}}+ n_{j}v^{j}\otimes t^{n_{j}-1}),
\end{eqnarray}
where $v^{j}\in U,\; n_{1}>\cdots >n_{s}$ and
$dv^{1}\otimes t^{n_{1}}+ n_{1}v^{1}\otimes t^{n_{1}-1}\ne 0.$
Compare the highest powers of $t$ on both sides of (\ref{ebasiscalculation}). 
If $dv^{1}\ne 0$, we have
$$u^{1}\otimes t^{m_{1}}=dv^{1}\otimes t^{n_{1}},$$
that is, $u^{1}=dv^{1}\in \im d$. This is a contradiction because
$0\ne u^{1}\in U'\cap \im d=0$.
Assume  $dv^{1}=0$. Then
$$u^{1}\otimes t^{m_{1}}=n_{1}v^{1}\otimes t^{n_{1}-1}$$
if $n_{2}\ne m_{1}-1$ and
$$u^{1}\otimes t^{m_{1}}=n_{1}v^{1}\otimes t^{n_{1}-1}+dv^{2}\otimes t^{n_{2}}$$
if $n_{2}= m_{1}-1$.
In both cases, we have $u^{1}\in U^{0}+\im d$,
noting that $v^{1}\in \ker d=U^{0}$. It is also a contradiction.

If $U^{0}\cap \im d=\ker d\cap \im d=0$, using a similar argument
we can prove that $\rho$ is a linear isomorphism from 
$(U^{0}\otimes t^{-1})\oplus (U'\otimes {\C}[t,t^{-1}])$ onto $L$.

(4) is obvious. 
\end{proof}

Now we turn our attention to a special class of local vertex Lie algebras
- -- local vertex Poisson algebras. 
First we recall the definition of Poisson algebras, which is needed in the
definition of vertex Poisson algebras.  

\bd{dpa}
A {\em Poisson algebra} is a commutative associative algebra $A$ equipped
with a Lie algebra structure on $A$ satisfying the following Leibniz rule:
\begin{eqnarray}\label{elebniz}
[a,bc]=[a,b]c+[a,c]b\;\;\;\;\mbox{ for }a,b,c\in A.
\end{eqnarray}
A {\em Poisson ideal} for a Poisson algebra $A$ is an ideal for both
the associative algebra $A$ and the Lie algebra $(A,[,])$.
\ed

Formula (\ref{elebniz}) is equivalent to that for $a\in V$, 
the operator $[a,\cdot]$
is a derivation of the associative algebra $A$.

Here is an example of Poisson algebras constructed from a Lie algebra.

\bex{exa0}
Let $\g$ be a Lie algebra. Then $S(\g)$, the symmetric algebra
over $\g$, could be considered as the algebra of polynomial functions
on $\g^{*}$, where $\g^{*}$ is considered as a manifold. 
It is well-known (cf. [W]) that $S(\g)$ has a Poisson algebra structure.
More specifically, let $u^{i}$ $(i\in I)$ be a basis of $\g$. Then we 
may identify $S(\g)$ with ${\C}[u^{i}\;|\;i\in I]$. For
$f,g\in {\C}[u^{i}\;|\;i\in I]$, we define
\begin{eqnarray}\label{epb0'}
\{f,g\}=\sum_{i,j\in I}\frac{\partial f}{\partial u^{i}}
\frac{\partial g}{\partial u^{j}}[u^{i},u^{j}].
\end{eqnarray}
Then $S(\g)$ becomes a Poisson algebra.
In general, for a manifold $M$, a Poisson algebra structure on $M$
(which means a Poisson algebra structure on $F(M)$, the space of 
smooth functions on $M$) amounts to a symplectic structure on $M$.
\eex

Let $(L,A,d,\rho)$ be a local vertex Lie algebra over  
a commutative associative algebra $A.$ 
Let $A(x)=\{a(x)\;|\;a\in A\}.$ Then $A\to A(x)$ via $a\mapsto a(x)$
is a linear isomorphism by Remark \ref{rdvla}. So $A(x)$
is a commutative associative algebra under product 
$$a(x)\circ b(x)=ab(x)$$
for $a,b\in A.$ Set
\begin{eqnarray}
R(x,y)=\sum_{k\in \N}A(x)\Delta^{(k)}(x,y)\subset A[[x,y,x^{-1},y^{-1}]].
\end{eqnarray}
Since $A$ has a unit $1$ by assumption, for $k\in \N$ we may identify
$\Delta^{(k)}(x,y)$ with $1\cdot\Delta^{(k)}(x,y)$, an element
of $R(x,y)$. Define an action of $A(x)$ on $R(x,y)$ by
\begin{eqnarray}
a(x)\circ \sum_{k\in \N}b_{k}(x)\Delta^{(k)}(x,y)
=\sum_{k\in \N}(ab_{k})(x)\Delta^{(k)}(x,y).
\end{eqnarray}
In view of Lemma \ref{lbasis}, the action is well defined and 
$R(x,y)$ is a free $A(x)$-module with a basis 
$\{\Delta^{(k)}(x,y)\;|\;k\in \N\}$.

\bd
(1)  A local vertex Lie algebra  $(L,A,d,\rho)$ is called 
a {\em local vertex Poisson algebra} if $A$ is a 
unital commutative associative algebra such that for $a,b\in A$,
\begin{eqnarray}
\{a(x),b(y)\}=\sum_{i}f_i(a,b)(y)\Delta^{(i)}(x,y),
\end{eqnarray}
where $f_{i}(a,b)\in A$ and that 
\begin{eqnarray}
&\{a(x),bc(y)\}=\{a(x),b(y)\}\circ c(y)+\{a(x),c(y)\}\circ b(y)\label{e3.7}
\end{eqnarray}
for $a,b,c\in A.$

(2) A local vertex Poisson algebra $(L,A,D,\rho)$ is called 
a local {\em vertex Poisson  differential algebra} 
if $D$ is a derivation of the associative algebra $A.$  
\ed

\br{bracket} $\{,\}$ in the definition of local vertex Poisson
algebra is called a vertex Poisson bracket on $A.$ In the following
when we talk about a vertex Poisson bracket on $A$ we always
mean a vertex Poisson bracket on $A$ associated to a local vertex
Poisson algebra over $A.$
\er

Here are some examples of local vertex Poisson algebras. 

\bex{ple1}
Let $A$ be a Poisson algebra. Set $L(A)=A\otimes\C[t,t^{-1}].$
Endow $L(A)$ with the loop Lie algebra structure of the Lie algebra
$(A,\{\cdot,\cdot\})$. From Example \ref{exa1} (3),
$(L(A), A, d, \rho)$ is a local vertex Lie algebra with
the domain of $d$ being 0, $\rho$ being the identity map and
$$\{a(x),b(y)\}=\{a,b\}(y)\Delta(x,y)\;\;\;\mbox{ for }a,b\in A.$$   
Let $a,b,c\in A$. We have
\begin{eqnarray}
& &\{ a(x),(bc)(y)\}=\{a,bc\}(y)\Delta(x,y)\nonumber\\
&=&(b\{a,c\})(y)\Delta(x,y)+
(\{a,b\}c)(y)\Delta(x,y)\nonumber\\
&=&b(y)\circ \{a,c\}(y)\Delta(x,y)+
(\{a,b\}(y)\Delta(x,y))\circ c(y)\nonumber\\
&=&b(y)\circ \{a(x),c(y)\}+\{a(x),b(y)\}\circ c(y).
\end{eqnarray}
Thus, $(L(A), A, d, \rho)$ is a local vertex Poisson Lie algebra.
\eex

\bex{exa2}
Let $\g$ be a Lie algebra and let $S(\g)$ be the 
symmetric algebra on $\g$.
Let $u^{i}$ $(i\in I)$ be a basis of $\g$. Then $S(\g)$ can be 
identified as the polynomial algebra ${\C}[u^{i}\;|\;i\in I]$ (see Example 
\ref{exa0}).  By Examples \ref{exa0}, \ref{exa1} (3) and \ref{ple1},
$L(S(\g))$ is a local vertex Poisson algebra over $S(\g).$ 
In particular, for $f,g\in {\C}[u^{i}\;|\;i\in I]$ we have 
\begin{eqnarray}\label{epb2}
 \{f(x),g(y)\}=\{f,g\}(y)y^{-1}\delta(\frac{x}{y})
=\left(\sum_{i,j\in I}\frac{\partial f}{\partial u^{i}}
\frac{\partial g}{\partial u^{j}}
[u^{i},u^{j}]\right)(y)y^{-1}\delta\left(\frac{x}{y}\right).\hspace{-1cm}
\end{eqnarray}
This vertex Poisson bracket on $S(\g)$ is usually called
the {\em ultra-Poisson} bracket.
\eex


\bp{p3.10}
Let $(L,A,D,\rho)$ be a local vertex Poisson differential 
algebra over $A.$ Then for any $a,b\in A,\;i\in {\N}$,
 there exists a  unique $a_{i}b\in A$ such that 
\begin{eqnarray}\label{d0}
[a(x),b(y)]=\sum_{i=0}^{\infty}\frac{1}{i!}(a_{i}b)(y)\Delta^{(i)}(x,y)
\end{eqnarray}
or equivalently
\begin{eqnarray}\label{d1}
[a(m),b(n)]=\sum_{i=0}^{\infty}\binom{m}{i}(a_{i}b)(m+n-i)\;\;\;\mbox{ for }m,n\in \Z.
\end{eqnarray}
Moreover $L^{0}=A/DA$ is a subalgebra of $L$
with $[a+DA,b+DA]=a_{0}b+DA$ for $a,b\in A$.
\ep

\begin{proof}
Note that $(Da)(x)=\frac{d}{dx}a(x)$ for any $a\in A$. Then in (\ref{edlva}) 
we can replace $f_{i}(u,v)^{(k_{i})}(y)$ by $D^{k_{i}}f_{i}(u,v)$.
For any $a,b\in A,\;i\in {\N}$, by Lemma \ref{lbasis}
we may define $a_{i}b\in A$ by requiring
\begin{eqnarray}
[a(x),b(y)]=\sum_{i=0}^{\infty}\frac{1}{i!}(a_{i}b)(y)\Delta^{(i)}(x,y).
\end{eqnarray}
Recall that $a(0)=\rho(a)$ and 
$L^0=\rho(A)$ which is isomorphic to $A/DA$ linearly. By  (\ref{d1}),
$[a(0),b(0)]=(a_0b)(0)$ for $a,b\in A$ or equivalently 
 $[a+DA,b+DA]=a_{0}b+DA.$ 
\end{proof}

Commutation relations (\ref{d0}) and (\ref{d1}) are exactly the commutator
relations for vertex operators in the theory of vertex operator algebras
(see [B], [FLM]). 

Next we discuss the relation between local vertex Poisson differential 
algebras and Poisson algebras. In the definition of local vertex Poisson 
differential algebra $A$ is only assumed to be a commutative associative 
algebra. From Proposition \ref{p3.10}, $A/DA$ has a Lie algebra 
structure already. If $DA$ is also an ideal of $A$ then $A/DA$
will be a Poisson algebra. But this is not true. It turns out that 
one has to modulo $A(DA)$ to get a Poisson algebra. 

\bp{pvpa}
Let $(L,A,D,\rho)$ be a vertex Poisson differential algebra. Then the quotient 
space $A/(DA)A$ is
a Poisson algebra with the original associative  multiplication and the
following Lie bracket:
\begin{eqnarray}
[a+(DA)A, b+(DA)A]=a_{0}b+(DA)A\;\;\;\mbox{ for } a,b\in A.
\end{eqnarray}
\ep

\begin{proof}
 It is clear that $(DA)A$ is an two-sided ideal of $A$, so that
$A/(DA)A$ is a commutative associative algebra. 
By Proposition \ref{p3.10},
$A/D(A)$ is a Lie algebra with the Lie bracket:
$[a+D(A),b+D(A)]=a_{0}b+D(A)$ for $a,b\in A$.
Taking ${\rm Res}_{x}{\rm Res}_{y}$ of (\ref{e3.7}) then gives
\begin{eqnarray}
[ab, c]=a[b,c]+b[a,c]\;\;\;\mbox{  for }a,b,c\in A.
\end{eqnarray}
Then $(DA)A/DA$ is a Lie ideal of $A/DA$ so that $A/(DA)A$ is a Lie algebra
with the Lie bracket, and the Leibniz rule is clear. The proof is 
complete. 
\end{proof}

Let $A$ be a polynomial algebra in variables $u_{i}^{(j)}$ for 
$i\in I,\;j\in {\N}$, where $I$ is an index set. Then $A$ has
a derivation uniquely determined by $\partial u_{i}^{(j)}=u_{i}^{(j+1)}$.
It is clear that any vertex Poisson bracket on $A$ is uniquely determined by
 $\{u_{i}(x),u_{j}(y)\}$ for $i,j\in I$.

\bd{ddgpb} 
A {\em differential-geometric} Poisson bracket on $A$ (see [DB]) is a 
vertex Poisson bracket $\{,\}$ on $A$ such that
\begin{eqnarray}
\{ u_{i}(x),u_{j}(y)\}=\sum_{l=0}^{n(i,j)}\psi_{ij}^{l}(u)(x)\Delta^{(l)}(x,y)
\end{eqnarray}
for any $i,j\in I,\; n(i,j)\in {\N}$, where $\psi_{ij}^{l}(u)\in A$.
A Poisson bracket on $A$  is said to be
{\em hydrodynamic type} if it has the form
\begin{eqnarray}\label{ehyd}
\{u^{i}(x),u^{j}(y)\}=g^{ij}(u(x))\Delta^{(1)}(x,y)
+\sum_{k\in I}b^{ij}_{k}(u(x))u^{(1)}_{k}
\Delta(x,y).
\end{eqnarray}
\ed

\br{rlinear}
In Definition \ref{ddgpb}, if each $\psi_{ij}^{k}(u)$ is a linear 
function in $u_{j}^{(l)}$ for $j\in I,\; l\in {\N}$, then one obtains
a (smaller) local vertex Lie algebra
over the space $U$ with a basis $\{u_{i}\;|\;i\in I\}$.
\er

The following interesting theorem was due to [DN]:

\bt{tDN}
(1) Under local changes of the fields $u=u(w)$ the coefficient $g^{ij}$ 
in the bracket (\ref{ehyd}) transforms like a bilinear form (a tensor 
with upper indices); if det $g^{ij}\ne 0$, then the expression 
$b_{k}^{ij}=g^{is}\Gamma_{sk}^{j}$ transforms in such a way that 
the $\Gamma_{sk}^{j}$ are the Christoffel symbols of a 
differential-geometric connection.

(2) In order that the bracket (\ref{ehyd}) be skew-symmetric it is necessary
and sufficient that the tensor $g^{ij}(u)$ be symmetric (i.e., that it 
defines a pseudo-Riemannian metric if det $g^{ij}\ne 0$) and the 
connection $\Gamma_{sk}^{j}$ be consistent with the metric, 
$g^{ij}_{k}=\Delta_{k}g^{ij}=0$.

(3) In order that the bracket (\ref{ehyd}) satisfy the Jacobi identity it
is necessary and sufficient that the connection $\Gamma_{sk}^{j}$ has no
torsion and the curvature tensor vanish. In this case the connection is 
defined by the metric $g^{ij}(u)$ which can be reduced to constant form.
\et

We discuss more examples:

\bex{exa5}
Suppose that $\{u_{i}(x),u_{j}(y)\}=d_{ij}\Delta^{(1)}(x,y)$ 
for $i,j\in I$, where 
$d_{ij}\in {\C}$. Set $H=\oplus_{i\in I} {\C}u_{i}$. Then
define $(u_{i}|u_{j})=d_{ij}$ for $i,j\in I$. 
Set $\hat{H}=H\otimes {\C}[t,t^{-1}]\oplus {\C}c$. Then a vertex Poisson 
bracket of this type determines a Lie algebra structure on $\hat{H}$
with $[u_{i}(m),u_{j}(n)]=d_{ij}m\delta_{m+n,0}c$ and $[\hat{H},c]=0$.
That is, a vertex Poisson bracket of this type determines a Heisenberg Lie algebra 
structure on $\hat{H}$.

More generally, suppose
\begin{eqnarray}
\{u_{i}(x),u_{j}(y)\}=d_{ij}\Delta^{(1)}(x,y)
+\sum_{k\in I}c_{ij}^{k}u_{k}(x)\Delta(x,y)
\end{eqnarray}
for $i,j\in I$, where $d_{ij}$ and $c_{ij}^{k}$ are constants.
Then a vertex Poisson algebra of this type gives rise to 
a Lie algebra structure on the subspace $\g$ with a basis 
$\{u_{1},u_{2},\ldots\}$ where
$$[u_{i},u_{j}]=\sum_{k\in I}c_{ij}^{k}u_{k}\;\;\;\mbox{ for }i,j\in I$$
and an affine Lie algebra structure on
$\hat{\g}={\g}\otimes {\C}[t,t^{-1}]+{\C}K$ with the standard 
bracket formula.
\eex

On the other hand, let $\g$ be a Lie algebra with
a nondegenerate symmetric invariant bilinear form $(,).$ Let 
$\{u_i\;|\;i\in I\}$ be a basis of ${\mathfrak g}$ such that
$(u_i,u_j)=d_{ij}$ and
$$[u_i,u_j]=\sum_{k}c_{ij}^ku_k.$$
Then the corresponding affine Lie algebra 
$\hat{\g}=\g\otimes {\C}[t,t^{-1}]+{\C}$
has bracket
$$\{u_{i}(x),u_{j}(y)\}=d_{ij}\Delta^{(1)}(x,y)
+\sum_{k\in I}c_{ij}^{k}u_{k}(x)\Delta(x,y)$$
where $u_i(x)=\sum_{n\in \Z}(u_i\otimes t^n)x^{-n-1}.$
Note also that $u(n)=u\otimes t^n$ in our notation.


\bex{eexa6} Consider a vertex Poisson Lie algebra structure on $A$
such that
\begin{eqnarray}\label{ee61}
\{u_{i}(x),u_{j}(y)\}=\sum_{k\in I}b^{k}_{ij}u_{k}^{(1)}(y)\Delta(x,y)
+\sum_{k\in I}c^{k}_{ij}u_{k}(y)\Delta^{(1)}(x,y).
\end{eqnarray}
If we define an algebra $B$ with a basis $\{a_{i}\;|\;i\in I\}$ such that
$a_{i}a_{j}=\sum_{k\in I}c_{ij}^{k}a_{k}$, then $B$ is commutative and 
associative [Za].
In this case, we have a local vertex Lie algebra over $B$ where
\begin{eqnarray}
[a(m),b(n)]=\frac{1}{2}(m-n)(ab)(m+n-1)\;\;\;\mbox{ for }a,b\in B,\;
 m,n\in {\Z}.
\end{eqnarray}
In [BN], essentially the same result was obtained by considering
\begin{eqnarray}\label{ee62}
\{u^{i}(x),u^{j}(y)\}
=\sum_{k\in I}d_{ij}^{k}u_{k}^{(1)}(x)\Delta(x,y)
+\sum_{k\in I}c_{ij}^{k}u_{k}(x)
\Delta^{(1)}(x,y),
\end{eqnarray}
where $d_{ij}^{k}, c_{ij}^{k}\in {\C}$ for $i,j,k\in I$.
Notice that using (\ref{echange})  we may exchange (\ref{ee62}) with 
(\ref{ee61}). The same result was also obtained in [GD].

One can also consider a vertex  Lie algebra over $B\oplus {\C}$ with
\begin{eqnarray}
\{u^{i}(x),u^{j}(y)\}
=\sum_{k\in I}d_{ij}^{k}u_{k}^{(1)}(y)\Delta(x,y)
+\sum_{k\in I}c_{ij}^{k}u_{k}(y)
\Delta^{(1)}(x,y)+f_{ij}\Delta^{(3)}(x,y),\nonumber\\
\end{eqnarray}
 where $d_{ij}^{k}, c_{ij}^{k}, f_{ij}\in {\C}$ for $i,j,k\in I$.
Let $B$ be defined as before and define a bilinear form $(\cdot|\cdot)$ on $B$
such that $(u_{i}|u_{j})=f_{ij}$ for $i,j\in I$.
Then
\begin{eqnarray}\label{elam}
[a(m),b(n)]=\frac{1}{2}(m-n)(ab)(m+n-1)+\frac{1}{6}(a|b)(m^{3}-m)\delta_{m,-n}c
\end{eqnarray}
for $a,b\in B,\; m,n\in {\Z}.$
It was  proved in [La] that (\ref{elam}) defines a Lie algebra structure
if and only if $B$ is commutative and associative, and
$(\cdot|\cdot)$ is a symmetric associative form on $B$.
(Lam obtained this result from a different point of view by considering
the weight 2 subspace of certain vertex operator algebras.)
\eex

\bex{eexa7} Consider a vertex Poisson algebra structure on $A$ such that
\begin{eqnarray}
\{u_{i}(x),u_{j}(y)\}=g^{ij}(x)\Delta^{(2)}(x,y)+
\sum_{k\in I}b^{ij}_{k}(x)u_{k}^{(1)}(x)\Delta^{(1)}(x,y).
\end{eqnarray}
Then it was proved in [Po] that the coefficients 
$g^{ij}$ and $b^{ij}_{k}$ must satisfy the following relations:
\begin{eqnarray}
& &g^{ij}=-g^{ji}, \;\; \frac{\partial g^{ij}}{\partial u_{k}}=b^{ij}_{k}, \;\;
 \sum_{l\in I}b^{ij}_{l}g^{lk}=\sum_{l\in I}b^{jk}_{l}g^{li};\\
& &\sum_{l\in I}\frac{\partial (b^{ij}_{l}g^{lk})}{\partial u_{m}}
=\sum_{l\in I}\left(b^{ij}_{l}b^{lk}_{m}
+b_{l}^{jk}b_{m}^{li}+b_{l}^{ki}b_{m}^{li}\right).
\end{eqnarray}

Furthermore, suppose that the coefficient $b_{k}^{ij}$ does not depend on 
$u_{i}$, then  $g^{ij}=\sum_{k\in I}b_{k}^{ij}u_{k}+g_{0}^{ij}$, 
where $b_{k}^{ij}$ are 
the structural constants of an associative anticommutative algebra $B$ that 
satisfy the condition
\begin{eqnarray}
\sum_{l\in I}b_{l}^{ij}g_{0}^{lk}=\sum_{l\in I}b_{l}^{jk}g_{0}^{li}.
\end{eqnarray}
Define a bilinear form on $B$ such that $(u_{i},u_{j})=g_{0}^{ij}$. Then
 we have
\begin{eqnarray}\label{eanti}
(u_{i}u_{j}, u_{k})=\sum_{l\in I}b_{l}^{ij}g_{0}^{lk}
=\sum_{l\in I}b_{l}^{jk}g_{0}^{li}=
(u_{j}u_{k},u_{i}).
\end{eqnarray}
Thus  $(ab,c)=(bc,a)$ for $a,b,c\in B$.
Using (\ref{echange}) we have
\begin{eqnarray}
\{u_{i}(x),u_{j}(y)\}
&=&\sum_{k\in I}\left(-b_{ij}^{k}u_{k}^{(1)}(y)\Delta^{(1)}(x,y)
+(b_{ij}^{k}u_{k}(y)+g^{ij}_{0})\Delta^{(2)}(x,y)\right)\nonumber\\
&=& -(u_{i}u_{j})^{(1)}(y)\Delta^{(1)}(x,y)
+ \left((u_{i}u_{j})(y)+(u_{i}|u_{j})\right)\Delta^{(2)}(x,y).
\end{eqnarray}
Thus
\begin{eqnarray}
& &[a(m),b(n)]\nonumber\\
&=&-m(m+n-1)(ab)(m+n-2)\nonumber\\
& &+m(m-1)\left((ab)(m+n-2)+(a|b)\delta_{m+n,1}\right))
\nonumber\\
&=&-mn(ab)(m+n-2)+m(m-1)\delta_{m+n,1}(a|b)
\end{eqnarray}
for $a,b\in B,\; m,n\in {\Z}$ and $[\hat{B}, {\C}]=0$.
One can easily show that $\hat{B}$ is a Lie algebra if and only if $B^{3}=0$.
\eex

\section{Constructing vertex algebras from a local vertex Lie algebra}

The main goal of this section is to
give a construction of vertex algebras from a local vertex 
Lie algebra. Roughly speaking, we shall construct
quantum objects, vertex algebras, from a classical object, 
a local vertex  Lie algebra. Results similar to some
of ours  have been also obtained in [K2], [MP], [Pr], [R] and [X].

First, let us review the definitions of vertex (operator) algebra
and module from [B] and [FLM] (cf. [Li1]).

\bd{dvoa}
A vertex algebra is a triple $(V,Y,{\bf 1})$ consisting of a 
vector space $V$,
a vector ${\bf 1}\in V$ and linear map $Y$ from $V$ to 
$({\rm End}V)[[x,x^{-1}]]$ satisfying the following axioms:

(1) $Y(u,x)v\in V((x))$ for $u,v\in V$,

(2) $Y({\bf 1},x)=1$, 

(3) $Y(u,x){\bf 1}\in V[[x]]$ and 
$\lim_{x\rightarrow 0}Y(u,x){\bf 1}=u$ for $u\in V$,

(4) The Jacobi identity holds for $u,v\in V$:
\begin{eqnarray}
& &x_{0}^{-1}\delta\left(\frac{x_{1}-x_{2}}{x_{0}}\right)Y(u,x_{1})Y(v,x_{2})
- -x_{0}^{-1}\delta\left(\frac{x_{2}-x_{1}}{-x_{0}}\right)Y(v,x_{2})Y(u,x_{1})
\nonumber\\
& &=x_{2}^{-1}\delta\left(\frac{x_{1}-x_{0}}{x_{2}}\right)Y(Y(u,x_{0})v,x_{2}).
\end{eqnarray}
\ed

We shall also use $V$ for the vertex algebra $(V,Y,{\bf 1})$.

\bd{dmodule}
Let $V$ be a vertex algebra. 
A $V$-module is a vector space $W$ equipped with a  linear map
$Y_W$ from $V$ to 
$({\rm End}\;W)[[x,x^{-1}]]$ satisfying the following axioms:

(1) $Y_W(u,x)w\in W((x))$ for $u\in V,\;w\in W,$

(2) $Y_W({\bf 1},x)=1,$ 

(3) The Jacobi identity holds for $u,v\in V$:
\begin{eqnarray}
& &x_{0}^{-1}\delta\left(\frac{x_{1}-x_{2}}{x_{0}}\right)
Y_W(u,x_{1})Y_W(v,x_{2})
- -x_{0}^{-1}\delta\left(\frac{x_{2}-x_{1}}{-x_{0}}\right)
Y_W(v,x_{2})Y_W(u,x_{1})
\nonumber\\
& &=x_{2}^{-1}\delta\left(\frac{x_{1}-x_{0}}{x_{2}}\right)
Y_W(Y(u,x_{0})v,x_{2}).
\end{eqnarray}
\ed

The notions of submodule  and homomorphism can be defined accordingly.

We shall need the following result from [Li1].

\bp{pli1}
Let $V$ be a vertex algebra, $W$ a $V$-module 
and $u$ a vector in $W$ such that $a_{n}u=0$ for any $a\in V,\; n\in {\N}$.
Then the linear map $\psi$ from $V$ to $W$ defined by $\psi(a)=a_{-1}u$ 
for $a\in V$ is a $V$-homomorphism.
\ep

A $V$-module $W$ is said to be {\em faithful} 
if the vertex operator map $Y_{W}$ is injective.
The following result was proved in [Li2].

\bp{px}
Let $V$ be a vertex algebra, $W$ a faithful $V$-module.
Let $a,b, c^{(0)},\ldots, c^{(k)}\in V$. Then 
$$[Y(a,x),Y(b,y)]=\sum_{j=0}^{k}\frac{(-1)^{j}}{j!}
Y(c^{(j)},y)\left(\dfrac{\partial}{\partial x}\right)^{j}
y^{-1}\delta\left(\frac{x}{y}\right)$$
if and only if
$$[Y_{W}(a,x),Y_{W}(b,y)]=\sum_{j=0}^{k}\frac{(-1)^{j}}{j!}
Y_{W}(c^{(j)},y)\left(\dfrac{\partial}{\partial x}\right)^{j}
y^{-1}\delta\left(\frac{x}{y}\right).$$
\ep

In [Li2], a general machinery has been built to produce vertex algebras
by using the notion of so-called local system of vertex operators.
Next we shall recall some results about local systems.

Let $M$ be a vector space. A {\it vertex operator} on $M$ is a formal series
$a(z)=\sum_{n\in {\Z}}a_{n}z^{-n-1}\in ({\rm End}\;M)[[z,z^{-1}]]$ such that
for any $u\in M$, $a_{n}u=0$ for sufficiently large $n$. All vertex 
operators on $M$ form a vector space (over ${\C}$), denoted by $VO(M)$.
On $VO(M)$, we have a linear endomorphism $D=\dfrac{d}{dx}$, the formal
differentiation.

Two vertex operators $a(z)$ and $b(z)$ on $M$ are said to be 
{\it mutually local} if there is a 
non-negative integer $k$ such that 
\begin{eqnarray}
(z_{1}-z_{2})^{k}a(z_{1})b(z_{2})=(z_{1}-z_{2})^{k}b(z_{2})a(z_{1}).
\end{eqnarray}
A space $S$ of vertex operators is said to be {\it local} if any two 
vertex operators of $S$ are mutually local, and a maximal local space of 
vertex operators is called a {\it local system}.

Let $V$ be a local system on $M$. Then $V$ is closed under the formal 
differentiation $D=\dfrac{d}{dx}$. For $a(x),b(x)\in VO(M)$, we define
\begin{eqnarray}\label{evp}
& &Y(a(x),z)b(x)\nonumber\\
&=&\Res_{z_{1}}\left(z^{-1}\delta\left(\frac{z_{1}-x}{z}\right)
a(z_{1})b(x)-z^{-1}\delta\left(\frac{x-z_{1}}{-z}\right)b(x)a(z_{1})\right).
\end{eqnarray}
Denote by $I(x)$ the identity endomorphism of $M$.

\bt{tli}
Let $M$ be a vector space and $V$ a local system on $M$. Then
$(V, Y, I(x))$ is a vertex algebra with $M$ as a natural module
such that $Y_M(a(x),z)=a(z)$ for $a(x)\in V.$
\et

Let $A$ be any local space of vertex operators on $M$. Then there exists a
local system $V$ that contains $A$. Let $\<A\>$ be the vertex subalgebra 
of $V$ generated by $A$. Since the vertex operator ``product'' (\ref{evp})
does not depend on the choice of local system $V$, $\<A\>$ is canonical.
Then we have:

\bc{clocal}
Let $M$ be a vector space and $A$ any local space of vertex 
operators on $M$. Then $A$ generates a canonical vertex algebra $\<A\>$
with $M$ as a natural module such that $Y_M(a(x),z)=a(z)$ for
$a(x)\in A.$
\ec

Next we should construct vertex algebras from any local 
vertex Lie algebra $L$.

\bp{padjoint}
Let $L$ be a local vertex Lie algebra over $U$, let $M$ be 
a restricted $L$-module and let $V$ be any vertex algebra of vertex operators 
on $M$ containing all $u(x)$ for $u\in U$. Then $V$ is an $L$-module
with $u(n)$ acting on $V$ by $u(x)_{n}$ for $u\in U,\; n\in {\Z}$.
\ep

\begin{proof} 
For $u,v\in U$, suppose 
$$[u(x),v(y)]=\sum_{j=0}^{n}\Delta^{(j)}(x,y)w_{j}(y),$$ 
where $w_{j}(y)$ lies in the space spanned by 
$w^{(i)}(y)$ for $w\in U$ and $i\geq 0.$ 
Then $V$ is an $L$-module if and only if 
\begin{eqnarray}\label{etep}
[Y(u(x),z_{1}),Y(v(x),z_{2})]
=\sum_{j=0}^{n}\Delta^{(j)}(z_{1},z_{2})Y(w_{j}(x),z_{2}).
\end{eqnarray}
Since $M$ is an $L$-module, we have
$$[\bar{u}(x),\bar{v}(y)]=\sum_{j=0}^{n}\Delta^{(j)}(x,y)
\bar{w_{j}}(y),$$
where $\bar{u}(x)=\sum_{m\in {\Z}}\bar{u}(n)$ and $\bar{u}(n)$ is an 
element of $\End M$ representing $u(n)$.
Since $M$ is a faithful $V$-module, (\ref{etep}) follows from 
Proposition \ref{px}. This concludes the proof. 
\end{proof}

Let $L$  be a local vertex Lie algebra over $U$. Recall from Proposition
\ref{central} (3) that $L$ has a polar decomposition $L=L^{+}\oplus
L^{-}$ where $L^{\pm}$ are Lie subalgebras of $L.$ Consider the
induced $L$-module $V(L)=U(L)\otimes_{U(L^+)}\C$ where $\C$ is the 
one-dimensional trivial $L^+$-module. Set ${\bf 1}=1\otimes 1.$
Then we have ([Li1], see also [K2], [FKRW], [MP], [Pr], [R], [X]):

\bt{tltov}
Let $L$ be any local vertex Lie algebra over $U$. Then
$V(L)$ has a vertex algebra structure with any restricted $L$-module
as a natural module.
\et

\begin{proof}
Let $M$ be any restricted $L$-module. Then $W=V(L)+M$
is a restricted $L$-module. Set $\bar{U}=\{ u(x)\;|\;u\in U\}$. Then
$\bar{U}$ is a local subspace of vertex operators on $W$. Thus $\bar{U}$ 
generates a canonical vertex algebra $V$ such that $W$ is a natural 
$V$-module. In particular, both $V(L)$ and $M$ are $V$-modules.
If we can prove that $V\simeq V(L)$, then we are done.
It is clear that $V$ is an $L$-module generated by $L$ from $Id_V$
by Proposition \ref{padjoint}. From the axioms of vertex
algebra we know that $u(n){\bf 1}=u(x)_n{\bf 1}=0$ for all $u\in U$
and $n\geq 0.$ Thus $V$ is a quotient $L$-module of $V(L)$.

On the other hand, $u(x)_n{\bf 1}=u(n){\bf 1}=0$ for all $u\in U$ and
$n\geq 0.$ It is easy to prove (cf. [GL]) that $a_n{\bf 1}=0$ for
all $v\in V$ and $n\geq 0.$ 
Then by Proposition \ref{pli1}, the linear map $\psi$ 
from $V$ to $V(L)$ defined by $f(a)=a_{-1}I(x)$ 
is a $V$-homomorphism. In particular
$\psi$ is an $L$-homomorphism. Thus $\psi$ 
is an $L$-isomorphism. The proof is complete. 
\end{proof}

Let $L$ be a local vertex Lie algebra over $U.$
Recall from Proposition \ref{padjoint} that $U^0$ is the kernal
of $d$ and $U^0(-1)=\{u(-1)\;|\;u\in U^0\}$ 
is a central subalgebra of $L.$ 
Let $\lambda$ be any linear character of $U^{0}(-1)$. 
Then we define $V(L,\lambda)$ to be the quotient module of $V(L)$ 
modulo the relation $c\cdot 1-\lambda(c)$ for $c\in U^{0}(-1)$.
As a simple consequence we have:

\bc{cltov}
Let $L$ be a local vertex Lie algebra over $U$ and $\lambda$
a linear character of $U^{0}$. Then $V(L,\lambda)$ is a quotient vertex
algebra of $V(L)$. 
\ec

Suppose that $L$ is a graded local vertex Lie algebra over
$U=\oplus_{n=0}^{\infty}U(n)$ such that 
$U(0)=U^{0}$. The degree of $u\in U(n)$ is defined to be $n.$ 
Then for any $\lambda\in U(0)^{*},$ 
$V(L,\lambda)$ has a natural ${\N}$-grading
$U(L,\lambda)=\oplus_{n\in {\N}}V(L,\lambda)_{(n)}$ with 
$$\deg(u_1(-n_1)\cdots u_k(-n_k){\bf 1}
=\deg\; u_1+n_1-1+\cdots +\deg\; u_k+n_k-1$$
for homogeneous $u^i$ in $U$ of positive degrees and with 
$V(L,\lambda)_{(0)}={\C}{\bf 1}$.
The $V(L,\lambda)$ has a unique maximal graded proper ideal $J(L,\lambda)$ 
so that $L(L,\lambda)=V(L,\lambda)/J(L,\lambda)$
is a simple graded vertex algebra. If the Virasoro algebra is a Lie subalgebra 
of $L$ associated with a vector $\omega\in U(2)$ and the ${\N}$-grading 
on $V(L,\lambda)$ is given by the operator $L(0),$ then $L(L,\lambda)$ 
is a simple vertex operator algebra.

\section{Lie algebras and Poisson algebras associated with
vertex operator algebras}

In this section, we study Poisson algebras and Poisson vertex algebras
associated with vertex algebras.

Let $(V,Y,{\bf 1})$ be a vertex algebra and let $D$ be the 
endomorphism of $V$ defined by 
\begin{eqnarray}
D(a)=a_{-2}{\bf 1}\;\;\;\mbox{ for }a\in V.
\end{eqnarray}
Set 
$\hat{V}=V\otimes {\C}[t,t^{-1}]$ and 
$\hat{D}=D\otimes 1+1\otimes \dfrac{d}{dt}$. 
Then $g(V)=\hat{V}/\hat{D}\hat{V}$ (by definition)
is a Lie algebra (cf. [B], [FFR], [Li3]) with bracket
\begin{eqnarray}\label{eb}
[a(m),b(n)]=\sum_{i=0}^{\infty}\binom{m}{i}(a_{i}b)(m+n-i),
\end{eqnarray}
where $a(m)$ stands for $a\otimes t^{m}+\hat{D}\hat{V}$.

Notice that in terms of generating functions, 
(\ref{eb}) can be written as
\begin{eqnarray}
[a(x),b(y)]
=\sum_{k=0}^{\infty}\frac{(-1)^{k}}{k!}(a_{k}b)(y)
\left(\dfrac{\partial}{\partial x}\right)^{k}y^{-1}\delta(\frac{x}{y})
\end{eqnarray}
where $a(x)=\sum_{n\in\Z}a(n)x^{-n-1}.$ 

We should also mention that any $V$-module becomes a natural $g(V)$-module 
with $a(m)$ acting as $a_{m}$ for $a\in V,\; m\in {\Z}$.

\br{rval}
Let $V$ be a vertex algebra and let $D$ be the endomorphism of $V$ 
defined above. Define an onto linear map $\rho$ from $\hat{V}$ 
to $g(V)$ by sending $u\otimes t^{n}$ to $u(n)$ for $u\in V,\; n\in {\Z}$. 
Then $\ker \rho=\im\hat{D}$. 
Thus the Lie algebra $g(V)$ is a local vertex Lie algebra over $V$.
\er

Let $V$ be a vertex algebra. Set 
\begin{eqnarray}
\hat{V}^{-}=V\otimes t^{-1}{\C}[t^{-1}],\;\;\;
\hat{V}^{+}=V\otimes {\C}[t],\;\;\;
\hat{V}^{0}=V\otimes {\C}\;(=V). 
\end{eqnarray}
Then
$\hat{V}^{\pm}$ and $V$ are $\hat{D}$-invariant subspaces.
Set $g(V)^{\pm}=\hat{V}^{\pm}/\hat{D}\hat{V}^{\pm}$ and $g(V)^{0}=V/DV$. Then 
it is clear that $g(V)^{\pm}$ and $g(V)^{0}$ are subalgebras of $g(V)$ and
we have the following polar decomposition 
\begin{eqnarray}
g(V)=g(V)^{+}\oplus g(V)^{-}.
\end{eqnarray}
Notice that ${\bf 1}(n)=0$ for $n\ne -1$ and that 
${\bf 1}(-1)$ is a central element of $g(V)$. Then for any 
$g(V)$-module $M$ 
and any complex number $\ell$, $({\bf 1}(-1)-\ell)M$ is a submodule of $M$.

Let $\ell$ be any complex number. Let ${\C}$ be the one-dimensional trivial 
$(g(V)^{+}\oplus g(V)^{0})$-module and consider the
induced module 
\begin{eqnarray}
I(V)=U(g(V))\otimes _{(g(V)^{+}\oplus g(V)^{0})}{\C}.
\end{eqnarray}
Set
\begin{eqnarray}
V^{[\ell]}=I(V)/({\bf 1}(-1)-\ell)I(V).
\end{eqnarray}
Then $V^{[\ell]}$ is a $g(V)$-module with ${\bf 1}(-1)$ 
acting as a scalar $\ell$. Let $\pi_{\ell}$ be the linear map
from $V$ to $V^{[\ell]}$ defined by
$\pi(a)=a(-1)1$ for $a\in V$ where $1=1\otimes 1+({\bf 1}(-1)-l)I(V).$
(We use the bold faced letter ${\bf 1}$ for the vacuum vector of $V.$)
Then it follows from the PBW theorem that $\pi_{\ell}$ is injective if 
$\ell\ne 0$.
If $\ell=1$, it follows from the universal property of $V^{[1]}$ that 
$V$ is a quotient $g(V)$-module of $V^{[1]}$.

\bt{tgv}
Let $V$ be any vertex algebra and let $\ell$ be any complex number. 
Then $V^{[\ell]}$ has a natural vertex algebra structure.
\et

\begin{proof} 
Since $g(V)$ is a local vertex algebra on the base 
space $V$, $V(g(V))$ 
is a vertex algebra. Then it is clear that $V^{[\ell]}$ is a quotient 
vertex algebra of $V(g(V))$ for any complex number $\ell$.
\end{proof}

Let $V$ be a vertex algebra and let $n$ be a positive integer. 
Then (see [FHL]) 
$V^{\otimes n}$ has a natural vertex algebra structure. Then $V^{\otimes n}$
is a natural $g(V)$-module of level $n$. Let $V^{(n)}$ be the $g(V)$-submodule
of $V^{\otimes n}$ generated from the vacuum ${\bf 1}$. It is easy to see that
$V^{(n)}$ is a vertex subalgebra of $V^{\otimes n}$.
It follows from the universal property of 
$V^{[n]}$ as a $g(V)$-module that
$V^{(n)}$ is a quotient vertex algebra of $V^{[n]}$.

The following lemma shows that $g(V)^{-}$ as a vector space
is linearly isomorphic to $V$.  (This result was also independently
obtained in [Pr].) 

\bl{llie}
Let $V$ be any vertex algebra. Then $V$ equipped with the product 
``$*$'' defined by
\begin{eqnarray}
a*b=a_{-1}b\;\;\;\mbox{ for }a,b\in V.
\end{eqnarray}
is a Lie admissible algebra with identity in the sense that
$V$ is a Lie algebra with the Lie bracket
\begin{eqnarray}
[a,b]=a*b-b*a=a_{-1}b-b_{-1}a\;\;\;\mbox{ for }a,b\in V.
\end{eqnarray}
Furthermore, it is isomorphic to $g(V)^{-}$.
\el

\begin{proof} 
Noticing that $(Da)(m)=-ma(m-1)$ for $a\in V,\;m\in {\Z}$, by using 
induction on $n$ we obtain
$a(-n)=\frac{1}{(n-1)!}(D^{n-1}a)(-1)$ for $a\in V$ and any positive 
integer $n$. Because of this relation, $g(V)^{-}$ is 
linearly spanned by $a(-1)$ for $a\in V$. Then the linear map $f$ from $V$ 
to $g(V)$ defined by $f(a)=a(-1)$ for $a\in V$ is onto.
Since $V$ is a $g(V)^{-}$-module with $a(-n)$ 
being represented by $a_{-n}$
for $a\in V,\; n=1,2,\ldots$, all $a_{-1}$ for $a\in V$ form
a Lie subalgebra of $gl(V)$.
Let $\pi$ be the representation on $V$ of $g(V)^{-}$.
Then $\pi f$ is the linear map from $V$ to 
${\rm End}V$ such that
$\pi f(a)=a_{-1}$. Because $a_{-1}{\bf 1}=a$ for any $a\in V$, $\pi f$ is 
injective, so that $f$ is injective. Thus $f$ is a linear isomorphism.
For $a,b\in V$, suppose $[a_{-1},b_{-1}]=c_{-1}$ for some $c\in V$.
Then
\begin{eqnarray}
c=c_{-1}{\bf 1}=[a_{-1},b_{-1}]{\bf 1}=a_{-1}b_{-1}{\bf 1}-b_{-1}a_{-1}{\bf 1}
=a_{-1}b-b_{-1}a.
\end{eqnarray}
Therefore, $[a_{-1},b_{-1}]=(a_{-1}b-b_{-1}a)_{-1}$. Then the lemma follows
immediately.
\end{proof}

\br{rfaith}
It follows from the proof of Lemma \ref{llie} that $V$ is a faithful
$g(V)^{-}$-module.
\er

Suppose that $V=\oplus_{n\in {\Z}}V_{(n)}$ is  a vertex operator 
algebra. For $a\in V_{(n)}$, we define
$\deg \; a=\wt\; a=n$. For $a\in V_{(m)},\; b\in V_{(n)}$, we have 
(see [FLM] for example)
$$\wt(a_{-1}b)=\wt(b_{-1}a)=\wt\; a+\wt\; b=m+n.$$
Then $V$ ($\simeq g(V)^{-}$) becomes a ${\Z}$-graded 
Lie algebra, so that
$g(V)^{-}$ and $V$ are isomorphic as graded vector spaces.

If $V$ is of CFT type in the sense that 
$V_{(n)}=0$ for all $n<0$ and $V_{(0)}={\C}{\bf 1}$ (cf. [DLMM]), 
then it follows from the PBW theorem that
$V^{[\ell]}$ is an ${\N}$-graded vector space with finite-dimensional 
homogeneous subspaces.

\bc{cgv}
Let $V$ be a vertex operator algebra of CFT type of rank $r$ and
$\ell$ a complex number. 
Then $V^{[\ell]}$ is a vertex operator algebra
of rank $r\ell$.
\ec

\br{raut}
Let $G$ be an automorphism group of a vertex operator algebra $V$.
Then $G$ is also an automorphism group for both the graded nonassociative  
algebra $(V,*)$ and the graded Lie algebra $(V,[,])$.
\er

\bl{center}
Let $V$ be a simple vertex operator algebra. Then 
the center of the Lie algebra or the nonassociative algebra $V$ 
is ${\C}{\bf 1}$.
\el

\begin{proof} 
It was proved in [Li1] that if $V$ is a simple vertex operator 
algebra, then $\ker L(-1)={\C}{\bf 1}$. Since $(V,*)$ is a graded algebra, the 
center is also graded. Let $u\in V$ be any homogeneous vector in the center 
of $V$.
Then $u_{-1}v=v_{-1}u$ for any $v\in V$. By the skew-symmetry we get
\begin{eqnarray}
& &u_{-1}v={\rm Res}_{z}z^{-1}Y(u,z)v
={\rm Res}_{z}z^{-1}e^{zL(-1)}Y(v,-z)u\nonumber\\
& &=\sum_{i=0}^{\infty}\frac{1}{i!}L(-1)^{i}(-1)^{i}v_{i-1}u.
\end{eqnarray}
Thus $\sum_{i=1}^{\infty}\frac{1}{i!}(-1)^{i}L(-1)^{i}v_{i-1}u=0$ 
for any $v\in V$.
Since $\ker L(-1)={\C}{\bf 1}$, we obtain
$$\sum_{i=1}^{\infty}\frac{1}{i!}(-1)^{i}L(-1)^{i-1}v_{i-1}u\in V_{(0)}$$
for any $v\in V$, so that $\wt\; u+\wt\; v-1=0$ if $v_{j}u\ne 0$ for some 
$j\in {\N}$. That is, if $\wt\; u+\wt v\;\ne 1$, then  $v_{j}u=0$ for any
 $j\in {\N}$.
Let $k$ be a positive integer such that $k+1>-\wt\; u$, so that
$\wt L(-1)^{k}\omega+\wt\; u >1$. Then we have
$$(L(-1)^{k}\omega)_{j}u=0\;\;\;\mbox{ for any }j\in {\N}.$$ 
In particular, $(L(-1)^{k}\omega)_{k}u=0$. Thus $\omega_{0}u=L(-1)u=0$.
Then $u\in {\C}{\bf 1}$. This concludes the proof.
\end{proof}

Another subalgebra $g(V)^{0}$ of $g(V)$ by definition is isomorphic to
the quotient space $V/DV$. If $V$ is a vertex operator algebra,
then $g(V)^{0}$ is a ${\Z}$-graded Lie algebra with $\deg a(0)=\wt a-1$ 
for any homogeneous vector $a\in V$.

\br{rtriangle}
For Kac-Moody algebras (see [K1]), one has the standard triangular 
decomposition  
where the positive part and the negative part are isomorphic (graded) 
Lie algebras.
If $V$ is a vertex operator algebra, one can 
get a triangular decomposition 
$g(V)=g(V)_{-}\oplus g(V)_{0}\oplus g(V)_{+}$ by defining 
$\deg a(n)=\wt a-n-1$
for any homogeneous element $a\in V$ and $n\in {\Z}$. Then one can prove
that $g(V)_{-}$ and $g(V)_{+}$ are isomorphic Lie algebras by using some 
of the proof for the contragredient module in [FHL]. However, $g(V)_{-}$
is not a graded Lie algebra with finite-dimensional homogeneous subspaces.
An interesting question is:  can we make $g(V)$ a doubly graded Lie algebra
with finite-dimensional homogeneous subspaces?
\er

Let $V$ be a vertex algebra as before. Set $P_n(V)=V/C_{n}(V)$
for $n\geq 2$ where $C_{n}(V)$ is the subspace of $V$
spanned by $a_{-n}b$ for $a,b\in V$. In the case that 
$V$ is a vertex operator algebra
associated to a highest weight module for the Virasoro algebra
or an affine Kac-Moody algebra, $P_n(V)$ has been studied in
[FF2] and [FKLMM]. We have (see [Zh]):

\bp{pzhu}
The space $P_2(V)$ is a Poisson algebra with the following associative
 multiplication
$\cdot$ and Lie multiplication $[,]$ defined by
\begin{eqnarray}
& &(a+C_{2}(V))\cdot (b+C_{2}(V))=(a_{-1}b+C_{2}(V)),\\
& &[a+C_{2}(V),b+C_{2}(V)]=a_{0}b+C_{2}(V)
\end{eqnarray}
for $a,b\in V$.
\ep

It follows from Proposition \ref{pzhu} that $C_{2}(V)$ is 
a two-sided ideal of the Lie admissible algebra $(V,*)$ 
such that the quotient algebra is commutative and associative.

An obvious question is that to what extent $P_2(V)$ determines $V$.  For
example, if $V$ and $U$ are two nonisomorphic vertex algebras, are
$P_2(V)$ and $P_2(U)$ nonisomorphic also?  Or more loosely, how much
information can we get from $P_2(V)$ for the vertex algebra
$V$?  It will be very nice if one can answer these questions. 

Next we shall compute $P_2(V(L))$ for a local vertex Lie algebra
$L$. Note that $V(L)$ is spanned by $u^1(-n_1)\cdots u^k(-n_s){\bf 1}$
for $s\geq 0,$ $u^i\in U,$ $n_i>0.$ In view of Remark \ref{rdvla} (1), 
the linear map $u\mapsto u(-1)$ from $U$ to $L$ is injective.
Then the linear map $u\mapsto u(-1){\bf 1}$ from $U$ to $V(L)$ is injective.
Now, we identify $U$ as a subspace of $V(L)$ 
through the linear map $u\mapsto u(-1){\bf 1}$ for $u\in U$.

Let $B$ be the subspace of
$V(L)$ linearly spanned by $a(-2-n)V(L)$ for $a\in U,\;n\in {\N}$.
Note that $a(m)=a_{m}$ for our notations.
{}From the definition of $C_{2}$, we have $B\subset C_{2}$. 
In the following we shall show that in fact, $B=C_{2}$. 
First, it follows from (\ref{ecomcom}) that $b(-1)B\subset B$ for $b\in U$.
Suppose $u\in V(L)$ such that 
$u_{-2-n}V(L)\subset B$ for all $n\geq 0.$ Let $a\in U, \;k\in {\N}$.
Then for $n\in {\N}$ using the Jacobi identity for $V(L)$  we have
\begin{eqnarray*}
& &{\rm Res}_{z}z^{-n-2}Y(a(-k-1)u,z)\nonumber\\
&=&{\rm Res}_{z}z^{-n-2}
{\rm Res}_{z_{1}}(z_{1}-z)^{-k-1}Y(a,z_{1})Y(u,z)\nonumber\\
& &-{\rm Res}_{z}z^{-n-2}{\rm Res}_{z_{1}}(-z+z_{1})^{-k-1}
Y(u,z)Y(a,z_{1})\nonumber\\
&=&{\rm Res}_{z}z^{-n-2}
\sum_{i=0}^{\infty}\binom{-k-1}{i}\left((-z)^{i}a(-k-1-i)Y(u,z)-
(-z)^{-k-1-i}Y(u,z)a(i)\right)\nonumber\\
&=&\sum_{i=0}^{\infty}\binom{-k-1}{i}\left((-1)^{i}a(-k-1-i)u_{-n-2+i}-
(-1)^{-k-1-i}u_{-n-k-3-i}a(i)\right)\nonumber\\
&=&a(-k-1)u_{-n-2}-(-1)^{-k-1}u_{-n-k-3}a(0)\nonumber\\
& &+\sum_{i\ge 1}\binom{-k-1}{i}\left((-1)^{i}a(-k-1-i)u_{-n-2+i}-
(-1)^{-k-1-i}u_{-n-k-3-i}a(i)\right).
\end{eqnarray*}
Then it follows that $(a(-k-1)u)_{-n-2}V(L)\subset B$.
(For $i=0$, we are using the fact that $b(-1)B\subset B$ and 
$b(-n-2)V(L)\subset B$ for $b\in U,\; n\in \N$.)
Since $V(L)$ is linearly spanned by $u^1(-n_1)\cdots u^s(-n_s){\bf 1}$
for $s\geq 0,$ $u^i\in U$ and $n_i>0$, it follows from induction that
 $u_{-n-2}V(L)\subset U$ for every $u\in V(L)$ and $n\in \N$.
Then from the definition of $C_{2}$, we have $C_{2}(V(L))\subset B$.
Hence $C_{2}(V(L))=B$.

Therefore, $P_2(V(L))=V(L)/B$ is spanned by 
$$u^1(-1)\cdots u^k(-1){\bf 1}+C_2(V(L))\;\;\;\mbox{ for }s\geq 0,\;u^i\in U.$$
 Set 
\begin{eqnarray}\label{elie}
{\mathfrak g}=\{u+C_2(V(L))\;|\;u\in U\}.
\end{eqnarray}
Then ${\mathfrak g}$ is a Lie subalgebra of 
$P_2(V(L))$ because from the bracket formula (\ref{ecomcom}), 
for $u,v\in U,$ 
\begin{eqnarray*}
& &[u+C_2(V(L)),v+C_2(V(L))]=u_0v+ C_2(V(L))\\
& &\ \ \ \ \ \  =(u_0v)_{-1}{\bf 1}+C_2(V(L))\\
& &\ \ \ \ \ \ =[u(0),v(-1)]{\bf 1}+C_2(V(L))\\
& &\ \ \ \ \ \ =\sum_{i, l_i=0}\binom{-1-l_{i}}{k_{i}}(-1)^{k_{i}}
f_i(u,v)(-1-l_i-k_i){\bf 1}+C_2(V(L))\\
& &\ \ \ \ \ \ =\sum_{i, l_i=k_i=0}f_i(u,v)(-1){\bf 1}+C_2(V(L)).
\end{eqnarray*}

Recall Proposition \ref{central} (3). Then we may get a basis of $L^{-}$
from $u^{0}(-1), u'(-n)$ for $u^{0}\in (U^{0})',\;u'\in U',\; n\ge 1$ 
by choosing a basis of $(U^{0})'$ and a basis of $U'$. Then from 
PBW theorem and $B=C_{2}(V(L))$ we immediately have:

\bl{lc2liealgebra}
Let $U'\subset U$ be such that $U=(U^{0}+\im d)\oplus U'$. Then 
there is a subspace $(U^{0})'$ of $U^{0}$ such that the map 
\begin{eqnarray}
& &(U^{0})'\oplus U'\rightarrow \mathfrak g\subset P_{2}(V(L))\nonumber\\
& &u\mapsto u+C_{2}(V(L))
\end{eqnarray}
is a linear isomorphism. Furthermore, if $\ker d\cap \im d=0$, 
we may take $(U^{0})'=U^{0}$.
\el

Now we are in the position to prove the following:
\bp{ppoisson}
Let $L$ be a local vertex Lie algebra over $U$, let $V(L)$ be 
the vertex algebra associated with $L$ and 
let $\mathfrak g$ be the Lie algebra defined by (\ref{elie}).
Then $P_2(V(L))$ is isomorphic to $S({\mathfrak g})$ as a Poisson algebra.
\ep

\begin{proof} The identity map on $\mathfrak g$  induces an
onto homomorphism of Poisson algebras from $S({\mathfrak g})$
to $P_2(V(L))$ by the universal mapping property of $S({\mathfrak g}).$ Since
$V(L)$ is isomorphic to $U(L^-)$ as a vector space, it follows
immediately from PBW theorem and $B=C_{2}(V(L))$ 
that $S({\mathfrak g})$ is, in fact, isomorphic to $P_2(V(L)).$
\end{proof}

We continue our discussion on the Lie algebra ${\mathfrak g}.$ 
Recall from Section 3 that
$L^0=\{u(0)\;|\;u\in U\}\subset L.$ 
It follows from Proposition \ref{central} (3) and Lemma \ref{lc2liealgebra} that
there is a linear map $\sigma$ from ${\mathfrak g}$
onto $L^0$ by sending $u+C_2(V(L))$ to $u(0)$ for $u\in U$. 
Since 
$$[u(0),v(0)]=\sum_{i, l_i=k_i=0}f_i(u,v)(0)+C_2(V(L)),$$
$\sigma$ in fact is a Lie algebra homomorphism. In general,
$\sigma$ is not an isomorphism. For example, if $u\in U^0$ 
and $u\not\in \im d$ then $u \not\in C_2(V(L))$ and
$u(0)=0.$

We need the following lemma:
\bl{lquot}
Let $V$ be a vertex (operator) algebra and let $I$ be an ideal of $V$. Then
$\bar{I}=(I+C_{2}(V))/C_{2}(V)$ is a Poisson ideal of $P_2(V)$ and 
$P_{2}(V/I)=P_2(V)/\bar{I}$.
\el

\begin{proof} 
Clearly, $\bar{I}$ is a Poisson ideal of $P_2(V)$. 
Furthermore, by definition we have 
$$P_{2}(V/I)=(V/I)/C_{2}(V/I),\;\;\;C_{2}(V/I)=(C_{2}(V)+I)/I.$$
Then $P_{2}(V/I)=P(V)/\bar{I}$.
\end{proof}

\bc{add} Let $\l$ be a linear character of $U^0.$ Then the Poisson algebra 
$P_2(V(L,\l))$ is isomorphic to $S(\mathfrak g)/J$ where $J$ is the ideal
generated by $(u-\l(u))+C_2(V(L))$ for $u\in U^0.$ 
\ec

\begin{proof} $V(L,\l))$ is a quotient vertex algebra of $V(L)$ modulo
the ideal $I$ generated by $u(-1)-\l(u)$ for $u\in U^0.$ Then $P_2(V(L,\l))$
is isomorphic to $S(\mathfrak g)/\bar{I}$ by Lemma \ref{lquot}. Clearly,
$J=\bar{I}$ is generated by $u-\l(u)+C_2(V(L))$ for $u\in U^0.$ 
\end{proof}.

Recall that for a vertex algebra $V,$ $L=g(V)$ is a vertex Lie algebra
over $V$ and that $g(V)^0=V/DV.$ 

\bc{cpoisson}
Let $V$ be a vertex algebra and let $\ell$ be a complex number. Then
$P_{2}(V^{[\ell]})$ is isomorphic to the Poisson algebra $S(V/DV)/J$
where $J$ is the ideal of $S(V/DV)$ generated by ${\bf 1}-\ell+DV.$ 
\ec
\begin{proof} 

First we deal with the case $\ell=0.$ Then ${\mathfrak g}$ in this 
case is isomorphic to $V/DV$ and   $P_{2}(V^{[0]})$
is isomorphic to $S(V/DV)$ by Proposition \ref{ppoisson}.
For general $\ell,$ the result follows from Corollary \ref{add}
and the fact that $V^{[\ell]}$ is a quotient vertex algebra 
of $V^{[0]}$ modulo the ideal generated by ${\bf 1}(-1)-\ell.$ 
\end{proof}

Let $\g$ be any Lie algebra.
 Set $L(\g)={\g}\otimes {\C}[t,t^{-1}]$. Then 
$L({\g})$ is a local vertex Lie algebra with bracket formula
\begin{eqnarray}
[a(m),b(n)]=[a,b](m+n)\;\;\;\mbox{ for }a,b\in {\g}.
\end{eqnarray}
Then $V(L({\g}))=U(\hat{\g}_{-})$ is a vertex algebra. Note that in this
case the ${\mathfrak g}$ in Proposition \ref{ppoisson} is isomorphic
to the Lie algebra ${\mathfrak g}$ we begin with.  
It follows 
{}from Proposition \ref{ppoisson} that $P_{2}(V(L({\g})))$ is 
isomorphic to $S({\g})$ as a Poisson algebra. Then we have:

\bc{clp}
Let ${\g}$ be any Lie algebra. Then there exists a vertex algebra $V$
such that $P_2(V)$ is isomorphic to $S({\g})$ as a Poisson algebra.
\ec




Our next goal is to compute $P_2(M(\ell,0))$ and $P_2(L(\ell,0))$
for affine vertex operator algebras $M(\ell,0)$ and 
$L(\ell,0)$ defined below.
Let $\mathfrak{g}$ be a simple Lie algebra and 
$H$ be a fixed Cartan subalgebra.
Denote by $\theta$ the longest root of $\mathfrak{g}$. Let 
$\hat{\mathfrak{g}}=\mathfrak{g}\otimes {\C}[t,t^{-1}]\oplus {\C}c$
be the affine Lie algebra. Set 
$\hat{\mathfrak{g}}_{\pm}=\mathfrak{g}\otimes t^{\pm 1}{\C}[t^{\pm 1}]$. Then 
$\hat{\mathfrak{g}}
=\hat{\mathfrak{g}}_{+}\oplus \mathfrak{g}\oplus {\C}c\oplus
 \hat{\mathfrak{g}}_{-}$.
For a complex number number $\ell$, set 
\begin{eqnarray}
M(\ell,0)=U(\hat{\mathfrak{g}})\otimes_{U(\hat{\mathfrak{g}}_{+}
+\mathfrak{g}+\C c)}\C,
\end{eqnarray}
the generalized Verma $\hat{\mathfrak{g}}$-module of level $\ell$. Let
$L(\ell,0)$ be the irreducible quotient module of $M(\ell,0)$. 
Then  $M(\ell,0)$ is a vertex operator algebra if $\ell$ is not negative the
dual Coxeter number (cf. [FZ], [Li2]) and  $L(\ell,0)$ is a
simple vertex operator algebra (cf. [DL], [FZ], [Li2]). 

\bp{pavg} The Poisson algebra $P_2(M(\ell,0))$ is isomorphic to
$S(\mathfrak{g})$. If $\ell$ is a positive integer, then
$P_2(L(\ell,0))$ is isomorphic to the quotient 
Poisson algebra of $S(\mathfrak{g})$ modulo the Poisson ideal 
generated by $e_{\theta}^{\ell+1}$.
\ep

\begin{proof} From the proofs of  Proposition \ref{ppoisson} 
and Corollary \ref{clp} we see that $P_2(M(\ell,0))$ is isomorphic to
$S(\mathfrak{g}).$

{}From the structure of $L(\ell,0)$ (cf. [K1]) we have
$L(\ell,0)= M(\ell,0)/I$, where $I=U(\hat{\mathfrak{g}}_{-})U(\mathfrak{g})
e_{\theta}^{\ell+1}(-1){\bf 1}$ and ${\bf 1}$ is a fixed highest
weight vector. Note that $I$ is an ideal of the vertex operator algebra 
$M(\ell,0)$.
Then $P_{2}(L(\ell,0))=P_{2}(M(\ell,0))/\bar{I}$ by Lemma \ref{lquot} and 
$\bar{I}=(I+C_2(M(\ell,0)))/C_{2}(M(\ell,0))$ is a 
Poisson ideal of $S(\mathfrak{g})$.
Since the associative product is defined as $a_{-1}b=a(-1)b$,
the Lie product is defined as $a_{0}b=a(0)b$ and 
$a(-n-2)b\in C_{2}(L(\ell,0))$ 
for $a\in \mathfrak{g},\;n\in {\N},\; b\in L(\ell,0)$,
it is clear that $\bar{I}$ is a subset of the Poisson ideal generated 
by $e^{\ell+1}_{\theta}$ of $S(\mathfrak{g})$.
Thus $\bar{I}$ is the Poisson ideal generated by $e_{\theta}^{\ell+1}$.
This concludes the proof.
\end{proof}


Finally we compute $P_2(V_L)$ for lattice vertex operator algebra $V_L.$ 
Let $L$ be a positive definite even lattice of rank $\ell$. 
Set 
\begin{eqnarray}
C_{2}(L)=\{\alpha\in L\;|\;
\<\alpha-\beta,\beta\> \le 0\;\;\;\mbox{ for all }\beta\in L\}.
\end{eqnarray}
Then $|\<\alpha,\beta\>|\le \<\beta,\beta\>$ for every $\beta\in L$ if 
$\alpha\in C_{2}(L)$.
Let $\alpha_{1},\ldots, \alpha_{n}$ be a ${\Z}$-basis for $L$ and let 
$\lambda_{1},\ldots,\lambda_{\ell}$ be the dual basis for the dual lattice 
$L^{o}$.
Let $k$ be a positive integer such that $kL^{o}\subset L$
Then for $\alpha\in C_{2}(L)$, we have
$|\< \alpha,\lambda_{i}\>|\le k \<\lambda_{i},\lambda_{i}\>$ 
for $i=1,\ldots, \ell$.
This implies that $C_{2}(L)$ is finite.
It follows from the definition that $0\in C_{2}(L)$ and 
$-\alpha\in C_{2}(L)$ if
$\alpha\in C_{2}(L)$. It also follows from the definition that 
$k\alpha \notin C_{2}(L)$ for $0\ne \alpha\in L,\; k\ne 0,\pm 1$.
Furthermore, $C_{2}(L)$ spans $L$ over $\Z$. If not, 
let $\alpha$ be an element of $L-\Z C_{2}(L)$ with the smallest length.
Since $\alpha\notin C_{2}(L)$, there is a $\beta\in L$ such that 
$\<\alpha-\beta,\beta\>\ge 1$. This implies that 
$\alpha-\beta\ne 0, \; \beta\ne 0$. From
\begin{eqnarray}\label{eexpression}
\<\alpha,\alpha\>=\<\alpha-\beta,\alpha-\beta\>+2\<\alpha-\beta,\beta\>+
\<\beta,\beta\>,
\end{eqnarray}
we get $|\alpha|>|\beta|, |\alpha-\beta|$. From the choice of $\alpha$,
we must have $\beta,\alpha-\beta\in \Z C_{2}(L)$. Thus $\alpha\in \Z C_{2}(L)$.
This is a contradiction.

\br{rkl}
In [KL], a subset $\Phi(L)$ of $L$ similar to $C_{2}(L)$ was introduced and
it was proved to satisfy the same properties.
\er

If $\alpha\in L-C_{2}(L)$, from the above argument and (\ref{eexpression}) we 
have $\<\alpha,\alpha\> \ge 6.$
Therefore $L_{2}\cup L_{4}\subset C_{2}(L)$, 
where $L_{n}=\{\alpha\in L\;|\;\<\alpha,\alpha\>=n\}.$
One can easily see that for $L={\Z}\alpha$, we have 
$C_{2}(L)=\{0, \pm \alpha\}$.

We shall need the explicit construction of $V_L$ given in [FLM]
including the group $\hat L$ and notation $\iota.$
Let $e$ be a section from $L$ to $\hat{L}$ such that $e(0)=1$ 
and $\epsilon(\cdot,\cdot)$ the corresponding 2-cocycle. 
Set 
\begin{eqnarray}
{\bf h}={\C}\otimes _{\Z}L,\;\;\;\hat{\bf h}={\bf h}\otimes {\C}[t,t^{-1}]\oplus {\C}.
\end{eqnarray}
As a vector space,
\begin{eqnarray}
V_{L}=S(\hat{\bf h}^{-})\otimes {\C}_{\epsilon}[L],
\end{eqnarray}
where $\hat{\bf h}^{-}={\bf h}\otimes t^{-1}{\C}[t^{-1}]$ and ${\C}_{\epsilon}[L]$
is the $\epsilon$-twisted group algebra of $L$.

For short we simply write $e^{\alpha}$
for $\iota(\alpha)$ for $\alpha\in L.$ 
Then $e^{\alpha}e^{\beta}=\epsilon(\alpha,\beta)e^{\alpha+\beta}$ for
$\alpha,\beta\in L$.
For $\alpha,\beta\in L$, from [FLM] we have
\begin{eqnarray}\label{evertexaction}
Y(e^\alpha,z)e^{\beta}=
\epsilon(\alpha,\beta)z^{\langle \alpha,\beta\rangle}
\exp \left(\sum_{m=1}^{\infty}\frac{\alpha(-m)}{m}z^{m}\right)
e^{\alpha+\beta}.
\end{eqnarray}
Then 
\begin{eqnarray}
& &(e^{\alpha})_{j}e^{\beta}=0\;\;\;\mbox{ for }
j> -1-\langle \alpha,\beta\rangle\\
& &(e^{\alpha})_{-\langle \alpha,\beta\rangle -1}e^{\beta}
=\epsilon(\alpha,\beta)e^{\alpha+\beta}. \label{edecompose}
\end{eqnarray}
Assume $\alpha\notin C_{2}(L)$. Then 
$\alpha=\beta_{1}+\beta_{2}$ for some $\beta_{1},\beta_{2}\in L$ 
with $\<\beta_{1},\beta_{2}\>\ge 1$.
{}From (\ref{edecompose}) we have
$$\epsilon(\beta_{1},\beta_{2})e^{\alpha}
=(e^{\beta_1})_{-1-\<\beta_{1},\beta_{2}\>}e^{\beta_{2}}
\in C_{2}(V_{L}),$$
noting that $-1-\<\beta_{1},\beta_{2}\>\le -2$.
Therefore
\begin{eqnarray}
e^{\alpha}\in C_{2}(V_{L})\;\;\;\mbox{  for }\alpha\in L-C_{2}(L).
\end{eqnarray}
Since $\alpha(-m)V_{L}\subset C_{2}(V_{L})$ for $m\ge 2$, from 
(\ref{evertexaction}) we have
\begin{eqnarray}\label{e5.16}
Y(e^{\alpha},z)e^{\beta}
\equiv \epsilon(\alpha,\beta)z^{\langle \alpha,\beta\rangle}
e^{\alpha(-1)z}e^{\alpha+\beta}\;\;\;\mbox{mod  }C_{2}(V_{L}).
\end{eqnarray}

For $\alpha\in L$, set 
\begin{eqnarray}
X_{\alpha}=e^{\alpha}+C_{2}(V_{L}),\;\;\;
Z_{\alpha}=\alpha(-1)+C_{2}(V_{L})\in P_{2}(V_{L}).
\end{eqnarray}
Then 
\begin{eqnarray}\label{ezoperatoradd}
Z_{\alpha+\beta}=Z_{\alpha}+Z_{\beta}\;\;\mbox{ for }\alpha,\beta\in L.
\end{eqnarray}

\bl{ladd}
The following relations hold in $P_{2}(V_{L})$ for $\alpha,\beta\in C_{2}(L)$:
\begin{eqnarray}
& &X_{\alpha}X_{\beta}=0 \;\;\mbox{ if }\alpha+\beta\notin C_{2}(L);
\label{edefiningb}\\
& &X_{\alpha}X_{\beta}=\frac{1}{(-\<\alpha,\beta\>)!}
\epsilon(\alpha,\beta)Z_{\alpha}^{-\<\alpha,\beta\>}X_{\alpha+\beta}
\;\;\mbox{ if } \alpha+\beta\in C_{2}(L);\label{exalphaxbeta}\\
& &Z_{\alpha}^{1-\<\beta-\alpha,\alpha\>}X_{\beta}=0.\label{ezpowersebeta}
\end{eqnarray}
In particular,
\begin{eqnarray}
& &Z_{\alpha}^{1+\<\alpha,\alpha\>}=0;\label{zpowers0}\\
& &Z_{\alpha}X_{\alpha}=0.\label{edefininge}
\end{eqnarray}
\el

\begin{proof} If $\alpha+\beta\notin C_{2}(L)$, we have
$e^{\alpha+\beta}\in C_{2}(V_{L})$. Since 
$\alpha(-m)C_{2}(V_{L})\subset C_{2}(V_{L})$ 
for $m\ge 1$, from (\ref{e5.16}) we have
\begin{eqnarray}\label{eprepaare1}
Y(e^{\alpha},z)e^{\beta}\equiv 0\;\;\;\mbox{mod  }C_{2}(V_{L}).
\end{eqnarray}
Then the relation (\ref{edefiningb}) follows immediately.

If $\alpha+\beta\in C_{2}(L)$, we have
$$\<\alpha,\beta\>=\<(\alpha+\beta)-\beta,\beta\>\le 0.$$
Using (\ref{e5.16})  we get
\begin{eqnarray}
(e^{\alpha})_{-1}e^{\beta}\equiv 
\frac{1}{(-\<\alpha,\beta\>)!}
\epsilon(\alpha,\beta)\alpha(-1)^{-\<\alpha,\beta\>}e^{\alpha+\beta}
\;\;\;\mbox{mod  }C_{2}(V_{L}).
\end{eqnarray}
This proves (\ref{exalphaxbeta}).

Using (\ref{e5.16}) again we get
\begin{eqnarray*}
(e^{\alpha})_{-2}e^{\beta-\alpha}\equiv 
\frac{1}{(1-\<\alpha,\beta-\alpha\>)!}
\epsilon(\alpha,\beta-\alpha)\alpha(-1)^{1-\<\alpha,\beta-\alpha\>}e^{\beta}
\;\;\;\mbox{mod  }C_{2}(V_{L}).
\end{eqnarray*}
Since $e^{\alpha}_{-2}e^{\beta-\alpha}\in C_{2}(V_{L})$, from this
we immediately obtain (\ref{ezpowersebeta}).
Clearly, (\ref{zpowers0}) and (\ref{edefininge}) are special cases of 
(\ref{ezpowersebeta}) with $\beta=0,\alpha$.
\end{proof}

Let $P(L)$ be the commutative associative algebra with identity
generated by symbols
$X_{\alpha}, Z_{\alpha}, \alpha\in C_{2}(L)$ with 
defining relations (\ref{ezoperatoradd})-(\ref{edefininge}). 
(With $C_{2}(L)$ being finite, it is clear that $P(L)$ is 
finite-dimensional.)

\bp{plattpoi}
The commutative associative algebra $P_{2}(V_{L})$ is isomorphic to 
the commutative associative algebra $P(L)$.
\ep

\begin{proof} From the proof of Lemma \ref{ladd} we see that $C_{2}(V_{L})$
contains the following subspaces:
\begin{eqnarray}
& &\alpha(-m-2)V_{L}\;\;\; \mbox{ for }\alpha\in L,\; m\in {\N};
\label{e1}\\
& &S(\hat{\bf h}^{-})\otimes e^{\beta}\;\;\mbox{
 for }\beta\notin C_{2}(L);\label{e2}\\
& &S(\hat{\bf h}^{-})
\alpha(-1)^{1-\<\beta-\alpha,\alpha\>}\otimes e^{\beta}\label{e4}
\end{eqnarray}
for $\alpha\in L,\; \beta\in C_{2}(L)$.
(Notice that $\<\beta-\alpha,\alpha\>\le 0$ for $\beta\in C_{2}(L)$.)
Let $U$ be the space of $V_{L}$ spanned by the above subspaces.
Noticing that $C_{2}(L)$ spans $L$ over $\Z$, we see that
$P_{2}(V_{L})$ as an associative algebra is generated by
$X_{\alpha}, Z_{\alpha}$ for $\alpha\in C_{2}(L)$ with the relations 
(\ref{ezoperatoradd})-(\ref{edefininge}).

To show that there is no more relation, we prove that $U=C_{2}(V_{L})$.
First, notice that the following holds:
$\alpha(-k-1)U\subset U$ for $\alpha\in L,\; k\in {\N}$. 
Because of (\ref{e1}) 
(\ref{e5.16}) is true with mod $C_{2}(V_{L})$ being replaced by mod $U$.
Then because of (\ref{e4}) it is clear that 
$(e^{\alpha})_{-n-2}e^{\beta}\in U$ for any 
$\alpha,\beta\in L,\; n\in \N$. 

Second, since 
$[\alpha(-k), (e^{\beta})_{-n-2}]
=\<\alpha,\beta\>(e^{\beta})_{-k-n-2}$ for $\alpha,
\beta\in L,\; k,n\in {\N}$, using induction we can prove
$$(e^{\alpha})_{-n-2}\alpha_{i_{1}}(-k_{1})\cdots \alpha_{i_{r}}(-k_{r})
e^{\beta}\in U$$
for any $\alpha,\alpha_{i},\beta\in L,\; k_{1},\ldots, k_{r}\ge 1$.

Third, similar to the proof of Proposition \ref{ppoisson}
using induction again we can prove that
$$\left(\alpha_{j_{1}}(-s_{1})\cdots \alpha_{j_{t}}(-s_{t})e^{\alpha}
\right)_{-n-2}\alpha_{i_{1}}(-k_{1})\cdots \alpha_{i_{r}}(-k_{r})
e^{\beta}\in U$$
for any $\alpha,\beta\in L,\; k_{1},\ldots, k_{r}\ge 1$. 
Thus $u_{-n-2}v\in U$ for $u,v\in V_{L},\; n\ge 0$, so that 
$C_{2}(V_{L}) \subset U$. Therefore $U=C_{2}(V_{L})$. 

{}From this, there exists a subset
$$G\subset 
\<\alpha(-1)\;|\;\alpha\in C_{2}(L)\>\otimes 
\{e^{\alpha}\;|\;\alpha\in C_{2}(L)\}
\subset V_{L}$$
such that $G\cap C_{2}(V_{L})=\emptyset$ and 
$V_{L}=C_{2}(V_{L})\oplus {\C}G$, 
where for a set $S$, $\< S\>$ denotes the free abelian semigroup 
(with identity) generated by $S$. 
Then $G$ gives rise to a basis of $P_{2}(V_{L})$.

In view of Lemma \ref{ladd}, we have an
algebra homomorphism $\psi$ from $P(L)$ onto $P_{2}(V_{L})$ such that
\begin{eqnarray}
\psi(Z_{\alpha})=\alpha(-1)+C_{2}(V_{L}),\;\;\; 
\psi(X_{\alpha})=e^{\alpha}+C_{2}(V_{L})
\;\;\;\mbox{ for } \alpha\in C_{2}(L).\hspace{-0.5cm}
\end{eqnarray}
Let $\tilde{G}$ be the corresponding subset of $P(L)$ with $\alpha (-1)$ and $e^{\alpha}$
being replaced by $Z_{\alpha}$ and $X_{\alpha}$, respectively.
Since $\psi(\tilde{G})=G$ in $P_{2}(V_{L})$ and 
$G$ gives rise to a basis of $P_{2}(V_{L})$, 
$\tilde{G}$ is linearly independent. It follows from the relations
(\ref{ezoperatoradd})-(\ref{edefininge}) that $\tilde{G}$ linearly spans $P(L)$.
Then it follows immediately that $\psi$ is a linear isomorphism. Therefore, $P_{2}(V_{L})$ is
isomorphic to $P(L)$ as an algebra.
\end{proof}

\bc{clatpo}
$P_2(L)$ has a Poisson algebra structure such that for $\alpha,\beta\in C_{2}(L)$,
\begin{eqnarray}
& &\{Z_{\alpha},Z_{\beta}\}=0,\;\; 
\{Z_{\alpha}, X_{\beta}\}=\<\alpha,\beta\>X_{\beta},\\
& &\{X_{\alpha},X_{\beta}\}=0\;\;\mbox{ if }\<\alpha,\beta\>\ge 0;\\
& &\{X_{\alpha},X_{\beta}\}=\frac{1}{(-\<\alpha,\beta\>-1)!}
Z_{\alpha}^{-\<\alpha,\beta\>-1}\epsilon(\alpha,\beta)X_{\alpha+\beta}
\end{eqnarray}
if $\alpha+\beta\in C_{2}(L)$ and $\<\alpha,\beta\>\le -1$.
\ec

\begin{proof} It follows from a calculation for the Lie brackets
in $P_{2}(V_{L})$. 
\end{proof}

\br{rnonpo}
Let $L$ be a nondegenerate even lattice that is not positive-definite. 
That is, there is an $\alpha\in L$ such that $\<\alpha,\alpha\>=k<0$.
Then $(e^{\alpha})_{k-1}e^{-\alpha}\in 1+C_{2}(V_{L})$. Thus
$1\in C_{2}(V_{L})$. Therefore, $P_{2}(V_{L})=0$. 
\er

Let $L={\Z}\alpha$ be a one-dimensional lattice such that $|\alpha|=2k$, 
where $k$ is a fixed positive integer. Then $C_{2}(L)=\{\alpha,0, -\alpha\}$.
Let $B_{k}$ be the quotient algebra
 of the polynomial algebra ${\C}[X,Y,Z]$ modulo the relations 
$X^{2}=Y^{2}=XZ=YZ=0,\; XY=\frac{1}{(2k)!}Z^{2k}$. Then $A_{k}$ is 
the algebra over the curve $X^{2}=Y^{2}=XZ=YZ=0,\; XY=\frac{1}{(2k)!}Z^{2k}$.
Define $\{Z,X\}=2kX,\; \{Z,Y\}=-2kY, \;\{X,Y\}=\frac{1}{(2k-1)!}Z^{2k-1}$. 
As an immediate consequence we have:

\bc{clattice}
The above defined $B_{k}$ is a Poisson algebra and it is 
isomorphic to  $P_{2}(V_{L})$.
\ec

\end{document}